\documentclass[11pt]{article}

\usepackage[colorlinks=false,pdfborderstyle={/W 1}]{hyperref}

\def\arXiv#1#2{\href{http://arxiv.org/abs/#1}{{\tt arXiv:#1 [#2]}}} 
\def\arXivo#1{\href{http://arxiv.org/abs/#1}{{\tt [arXiv:#1]}}} 

\usepackage{amsmath}
\usepackage{amssymb}
\usepackage{amsthm}
\usepackage{amsfonts}
\usepackage{graphicx}
\usepackage{mathrsfs}

\usepackage[usenames,dvipsnames]{xcolor}

  \ifx\LabelFigloaded\MYundefined\relax
  \else
   \fi

  \def\LabelFigloaded{\relax}

  \chardef\LabelFigCatAt\the\catcode`\@
  \catcode`\@=11

 \let\LabelFigwlog@ld\wlog
 \def\wlog#1{\relax}

 \ifx\\\MYundefined@
    \let\\\relax
 \fi

  \def\ms@g{\immediate\write16}

 \def\N@wif{\csname newif\endcsname }
 \def\Temp@ {\N@wif\ifIN@}
 \ifx\INN@\MYundefined@
    \else \let\Temp@\relax
 \fi
 \Temp@

  \def\IN@{\expandafter\INN@\expandafter}
  \long\def\INN@0#1@#2@{\long\def\NI@##1#1##2##3\ENDNI@
    {\ifx\m@rker##2\IN@false\else\IN@true\fi}%
     \expandafter\NI@#2@@#1\m@rker\ENDNI@}
  \def\m@rker{\m@@rker}
 
  \newtoks\Initialtoks@  \newtoks\Terminaltoks@
  \def\SPLIT@{\expandafter\SPLITT@\expandafter}
  \def\SPLITT@0#1@#2@{\def\TTILPS@##1#1##2@{%
     \Initialtoks@{##1}\Terminaltoks@{##2}}\expandafter\TTILPS@#2@}

 \def\Shifted@@#1#2#3{\setbox0=\hbox{#3}%
   \raise -\dp0\vbox {\kern-#2%
       \hbox {\kern#1\unhbox0\kern-#1}%
           \kern#2}}

 \newcount\gridcount
 \newbox\auxGridbox@ \newbox\hGridbox@ \newbox\vGridbox@
 \newbox\Labelbox@ \newbox\auxLabelbox@
 \newbox\Coordinatebox@
 \newtoks\Labeltoks@
 \newdimen\Wdd@ \newdimen\Htt@
 \newdimen\Wddd@ \newdimen\Httt@
 
 \def\Wr@{\immediate\write16}

 \newdimen\GL@wd
 \def\GridLineWidth#1{\GL@wd=#1}

 \def\gobble#1{}
 \def\EdgeErr@{\Wr@{}%
      \Wr@{\string\Edges\space argument
      1, 10, 100 or 1000 please\string!}%
      }

 \newcount\Edgect@

 \def\Sweepup#1\endSweepup{}

 \def\SetEdges@{%
    \edef\Zr@@s{\expandafter\gobble\number\Edgect@\empty}%
        \count255=0\Zr@@s\relax
        \ifnum\count255=\z@\else\EdgeErr@\show\tailtest\fi
        \count255=1\Zr@@s\relax
        \ifnum\count255=\Edgect@\relax\else\EdgeErr@\show\leadtest\fi
    \EdgGl@b\edef\Zr@s{\expandafter\gobble\Zr@@s\empty}
    \ifnum\Edgect@>\@ne\relax\EdgGl@b\let\L@Dc\empty
        \else\EdgGl@b\edef\L@Dc{\string.}\fi
    \ifnum\Edgect@>\@ne\relax
        \EdgGl@b\edef\Edgescale@##1{\divide##1 by \Edgect@}%
        \else\EdgGl@b\edef\Edgescale@##1{}\fi
    }

 \def\Edges#1{\Edgect@=#1\relax
     \let\EdgGl@b\global \SetEdges@}

 \Edges{1}

 \def\hhrule{\hrule height \GL@wd\vskip-.\GL@wd}

 \def\hRule@{%
   \advance\gridcount -2%
   \vfil\hhrule\vfil
   \llap{\smash{\raise -2.5pt
     \hbox{\L@Dc\number\gridcount\Zr@s\kern2pt}}}%
   \hhrule
   }

\def\vvrule{\vrule width \GL@wd \kern-\GL@wd}

 \def\vRule@{\advance\gridcount 2%
   \hfil\vvrule\hfil
   \setbox\auxGridbox@=\vbox to 0pt
      {\vskip \Htt@\vskip 2pt
        \hbox to 0pt{\hss\L@Dc\number\gridcount\Zr@s\hss}\vss}%
      \wd\auxGridbox@=0pt \box\auxGridbox@
   \vvrule
   }

 \def\PlaceGrid@@{\gridcount=10 
  \setbox\hGridbox@=\hbox{%
        \hbox{%
             \hskip-.4pt\vrule
             \vbox to \Htt@{%
               \offinterlineskip\parindent=\z@\relax
               \hbox to \Wdd@{\hfil}
               \hRule@\hRule@\hRule@\hRule@
               \vfil\hhrule\vfil}%
             \vrule\hskip-.4pt}
    }%
  \gridcount=0%
  \setbox\vGridbox@=\hbox{%
      \vbox{\offinterlineskip\parindent=0pt\hsize=0pt
         \vskip-.4pt\hrule%
         \hbox to \Wdd@{%
                 \vtop to \Htt@{\vfil}%
                 \vRule@\vRule@\vRule@\vRule@
                 \hfil\vvrule\hfil}%
         \hrule\vskip-.4pt}}%
  \wd\hGridbox@=0pt\ht\hGridbox@=0pt
  \wd\vGridbox@=0pt\ht\vGridbox@=0pt
  \hbox{\box\hGridbox@\box\vGridbox@}%
  }

 \def\LabelsGlobal{\def\LabGl@b{\global}}
 \def\LabelsLocal{\def\LabGl@b{}}
 \LabelsGlobal 

 \def\SetLabels#1\endSetLabels{%
   \LabGl@b\Labeltoks@={#1()\\}%
   }

 \LabGl@b\Labeltoks@={()\\}

 \def\ShowGrid{\LabGl@b\let\PlaceGrid@\PlaceGrid@@}
 \def\HideGrid{\LabGl@b\let\PlaceGrid@\relax}
 \def\Grids{\ShowGrid\LabGl@b\let\GridSwitch@\ShowGrid}
 \def\noGrids{\HideGrid\LabGl@b\let\GridSwitch@\HideGrid}

 \noGrids

 \def\bAdjust@@{%
     \setbox\auxLabelbox@=\hbox{\raise \dp\auxLabelbox@
            \box\auxLabelbox@}}
 \def\bAdjust@{\let\vAdjust@\bAdjust@@}

 \def\eAdjust@@{\dimen0=-.5\ht\auxLabelbox@
     \advance\dimen0 by .5\dp\auxLabelbox@
     \setbox\auxLabelbox@=
            \hbox{\raise\dimen0\box\auxLabelbox@}}
 \def\eAdjust@{\let\vAdjust@\eAdjust@@}

 \def\tAdjust@@{%
     \setbox\auxLabelbox@=\hbox{\raise-\ht\auxLabelbox@
            \box\auxLabelbox@}}
 \def\tAdjust@{\let\vAdjust@\tAdjust@@}

 \let\vAdjust@\relax

 \def\lAdjust@{\let\hAdjust@\rlap}
 \def\rAdjust@{\let\hAdjust@\llap}

 \let\hAdjust@\relax\let\vAdjust@\relax

 \def\FetchLabel@#1(#2)#3\\{%
     \IN@0#2@@\ifIN@
        \setbox0=\hbox{\ignorespaces#1#3\unskip}%
        \ifdim\wd0>0pt
           \ms@g{}%
           \ms@g{ !!! Bad label(s)? !!!}%
           \message{ #1(#2)#3}%
        \fi
        \def\LabelMole@##1\endFetchLabel@{%
            \IN@0()\\@##1@%
            \ifIN@\def\Temp@{\FetchLabel@##1\endFetchLabel@}%
            \else\def\Temp@{}%
            \fi
            \Temp@
           }%
     \else
       \ignorespaces#1\unskip
       \setbox\auxLabelbox@=%
         \hbox to 0pt{\hss\ignorespaces\hAdjust@
          {\ignorespaces#3\unskip}\hss}%
       \vAdjust@
       \let\hAdjust@\relax\let\vAdjust@\relax
       \AugmentLabelBox@@{#2}%
       \ht\Labelbox@=0pt\dp\Labelbox@=0pt
       \let\LabelMole@\FetchLabel@%
     \fi\LabelMole@}

 \newtoks\XYSep@ 
 \def\SetXYSeparator#1{%
     \IN@0#1@@\ifIN@\XYSep@{*}%
     \else
     \XYSep@{#1}%
     \fi
     }

 \SetXYSeparator*

 \def\AugmentLabelBox@@#1{%
     \IN@0\the\XYSep@ @#1@\ifIN@
       \SPLIT@0\the\XYSep@ @#1@%
       \setbox\Labelbox@=\hbox to 0pt{%
         \unhbox\Labelbox@
         \Shifted@@{\the\Initialtoks@\Wddd@}%
         {\the\Terminaltoks@\Httt@}%
         {\box\auxLabelbox@}}%
     \else
         \ms@g{}%
         \ms@g{ !!! Bad insertion point. !!!}%
         \message{ (#1\ this point was rejected.)}%
     \fi
    }

 \def\FetchOption@#1[#2]#3\endFetchOption@{%
    \def\temp{#1}
    \ifx\temp\empty
       \Edgect@=#2\relax
       \let\EdgGl@b\relax
       \SetEdges@
       \Cleaner@#3%
    \fi}

 \def\Cleaner@#1[@]{\Labeltoks@{#1}}
     
 \def\PlaceLabels@@{\mathsurround=0pt
     \def\Cr@{\\}%
     \let\L\lAdjust@\let\R\rAdjust@
     \let\B\bAdjust@\let\E\eAdjust@\let\T\tAdjust@
     \expandafter\FetchOption@\the\Labeltoks@[@]\endFetchOption@
     \Wddd@=\Wdd@ \Edgescale@\Wddd@ 
     \Httt@=\Htt@ \Edgescale@\Httt@
     \expandafter\FetchLabel@\the\Labeltoks@\endFetchLabel@
     \box\Labelbox@
     }%

 \let \PlaceLabels@\PlaceLabels@@

 \def\AffixLabels#1{\setbox\Coordinatebox@=\hbox{#1}%
      \Wdd@=\wd\Coordinatebox@ \Htt@=\ht\Coordinatebox@
      \advance\Htt@ \dp\Coordinatebox@
      \hbox{\copy\Coordinatebox@\kern-\Wdd@ 
           \Shifted@@{0pt}{-\dp\Coordinatebox@}%
           {\PlaceLabels@\PlaceGrid@}%
           \kern\Wdd@}%
      \GridSwitch@ 
      \LabGl@b\Labeltoks@{()\\}%
      }
 
   \let\wlog\LabelFigwlog@ld   
   \catcode`\@=\LabelFigCatAt  

                                By

              Raymond S\'eroul <A18645@FRCCSC21.BITNET>
                                and 
              Laurent Siebenmann <lcs@topo.math.u-psud.fr>
    
              VERSIONS: July 1991, Oct 1991, Jan 1992, July 1992

INTRODUCTION

      This labelling package is intended for TeX users who
rely on non-TeX sources for for their graphics inserts.  It
provides means for adding TeX labels to such inserts with a
minimum of fuss. 

       For most labels, TeX users have in the past found it
reasonably convenient to rely on non-TeX sources. Typical
occasions when an inescapable need for TeX labels seemed to
arise are

 (a) when the graphics program lacks certain exotic or complex
mathematical symbols

 (b) when the very highest typographical quality is wanted for the
labels

 (c) when labels included with the graphics fail to print, 
 and you cannot figure out why (cf. boxedeps.doc).  The labels
 provided by labelfig.tex are 100

       Since this package first appeared, many users, who in the
past scarcely dreamed of using TeX labels, have come to use
nothing but.  So it is now appropriate to add

Intoxication Warning:  TeX labels may be addictive and expensive. 

     If you have a fast preview you may disagree, and even find
that this package provides an agreeable paste-up environment; see
extra applications at end.

     Note to publishers: It is possible and convenient to ultimately
export the TeX labels produced by labelfig.tex to become an integral
part of the EPS file. This is often desired by a publisher who typically
uses an "upmarket" graphics or page layout program, with which the
staff is skilled in perfecting figures.  See Appendix I for
a recipe.

     The authors are grateful to Patrick Ion of Math Reviews for
helpful comments and encouragement.

BASIC INSTRUCTIONS

    After reading in the macro file using

preview or proof your figure with a coordinate grid printed on
top, by typing the following:

    \ShowGrid  
    \AffixLabels{<the graphics insertion>}

Here <the graphics insertion> is what you would type to insert
the graphics object alone without the grid.  This must provide
for the space around it. For example <the graphics insertion>
might well be \BoxedEPSF{MyFigure scaled 700} using the
boxedeps.tex macro package (from same source); this provides a
TeX box containing the encapsulated PostScript insert specified by
the file MyFigure. \AffixLabels{...} provides the grid (supposing
\ShowGrid is present) and later, once you have specified labels
using the grid, it will "tack on" the labels.

     The grid is a sort of (usually elongated) checkerboard of
ten rows and ten columns and its (internal) partitions are by
default numbered  .1, ... ,.9  both horizontally (X-coordinate
running left to right) and vertically (Y-coordinate running bottom
to top).  Thus the points enclosed by the grid correspond to the
points of the unit square in the cartesian "X-Y" plane, the lower
left corner corresponding to the origin (0,0).  By extrapolation,
the full page corresponds to a larger rectangle in the plane.

     These coordinates serve to position labels as follows.
Before the \AffixLabels{...} command type label specifications:

  \SetLabels
   (<X-coordinate>*<Y-coordinate>) <first label> \\
   .
   .
   .
   (<X-coordinate>*<Y-coordinate>)  <last label> \\
  \endSetLabels

Each row specifies one label and is terminated by \\.  In each
row, the position indicator comes first; it is written as a
standard cartesian point except that the X- and Y- coordinates
are separated by * rather than a comma because TeX allows a
comma as decimal point. There are no dimension units to specify
as the unit is the grid itself.

     By default, this cartesian point specifies where the middle
of the baseline of the label will be located.  However if you precede
the point by \L [or \R] the left [or right] edge of the baseline will
be located there. Similarly you may also precede the point by \T, \E,
or \B to vertically align the top equator or bottom of the label box
at the specified point.  This gives nine standard positions of
the label with respect to the insertion point --- corresponding to
the eight principle points of the compas and the center

                     \L\T     \T      \R\T

                     \L\E     \E      \R\E

                     \L\B     \B      \R\B

But this neglects the default "baseline" level of TeX,
giving potentially three more positions

                     \L    <no tag>   \R

For text, the baseline level is often the preferred. Its relation to
the others is variable. It will often coincide with the bottom level,
as happens for "X".  But it is often distinct, as for "g", in which
case you have in all 12 distinct positions rather than 9.

     It is convenient to think of this specification of label
position as attaching the label by a thumb-tack to the coordinate
grid. There are up to twelve positions of the thumb-tack on the
label, while the position of the thumb-tack on the coordinate grid is
arbitrary.  Normally, one choses the position of the thumb-tack on
the label to be the one that is the closest to the item being
labeled.  There are good reasons for this "rule of thumb":

   (a)  It facilitates correct positioning at first try.

   (b)  If the scale of the figure must be altered after labels
have been affixed, the labels have a good chance of remaining well
positioned.

   (c)  The visible grid need not extend beyond the "bounding box"
for the figure, because the best preferred position is always
(at least almost) within the bounding box .

The second reason is particularly important. Indeed it often
happens that scale has to be altered after labelling begins, in
order to either provide space for the labels, or to adjust
proportions between the labels and the figure.  (The size of labels
is unaffected by scaling.)

     Here is an artificial but self-contained test which uses
TeX rules to make a graphics object.

TEST

    Do not skip this!

 \def\FrameIt#1{\hbox{\vrule$\vcenter {\hrule\kern3pt%
             \hbox {\kern3pt #1\kern3pt}%
               \kern3pt\hrule}$\relax\vrule}}

 \def\Caption#1#2{\FrameIt{%
       \vtop {\hsize=#1\relax \parindent=0pt
         \leftskip=0pt \rightskip=0pt plus15pt
         \parfillskip=0pt
         \lineskip=1pt\baselineskip=0pt
         #2}}}

 \def\FirstQuadrant{\hbox to 100pt{\vrule\vbox to 100pt{%
        \hbox to 100pt{\hfil}\vfil\hrule}\hss}}

  \SetLabels
    \R(.5*.2) $\zeta\,\cdot$\\
    (.9*-.10) $\xi$\\
    \R(-.03*.9) $\eta$\\
    \T(.5*.9) \Caption{70pt}{%
          \it The norm of
          $g(\xi+i\eta)$ is indicated on
          contours of this invisible surface.}\\
  \endSetLabels

  \AffixLabels{\FirstQuadrant} \end

  Note that the coordinates to use for labels are indicated on the
edges of the grid (when visible) corresponding to the conventional
x- and y- axes of the Cartesian plane. By default the grid is
1-by-1. However, by the command \Edges{100}, you can change this
to 100-by-100 and many users find this alternative most
convenient. Place the command \Edges{...} in your style file (or
header) since its effect is is global. Other possible edge values
are 10 and 1000.

  If you use the command \Edges{...} at all, do so with care.  For
if you accidentally delete an \Edges{...} command your labels will
abruptly be badly misplaced and may logically but mysteriously
generate "dimension too big" errors under TeX and "off page" errors
under your driver.  

  You can dictate the edgescale for an individual figure by giving
the scale in brackets immediately after \SetLabels.  Thus, to
import into an article using say \Edge{100} a figure labelled using
another edgescale, say the original 1-by-1 default, you can use
\SetLabels[1]...\endSetLabels.

GETTING IT DOWN PAT

     Complicated labeling deserves the same respect as
complicated mathematics.  Do not expect it to come out perfect the
first time!  What is needed in either case is a mechanism to
repeatedly typeset troublesome pieces.

     One mechanism is always available.  One does complicated
labelling in a separate "test" file involving just the figure being
labelled;  a texpert will know how to \dump TeX's current state as
a temporary format that restarts rapidly at each retry.  Usually,
one then pastes the completed labelled figure back into the main
TeX file, but, of course, one can also \input it as an auxiliary
file.

     If you do not have a TeXpert at handy, here is a first
approximation to an efficient setup. By deletions reduce a copy
of your article to just a few lines before and after the figure.
Now label the figure, and finally, copy and paste the labelled
figure to the original article. Then copy the next figure to label
into this testbed and repeat. The TeXpert can improve the  speed
at which TeX starts up, by compiling a format specifically for
your article; just one caution: best NOT include in the format
ephemeral details of setup like \Set<mydriver>ArtSpecials (from
boxedeps.tex because this reads  figure dimensions which you may
change during your work session.

     An improved mechanism to repeatedly typeset troublesome
pieces is now available on the Macintosh; it is called LinoTeX;
see the same ftp sources.  It could be set up on many types
of computer.

     Before using labelfig.tex to attach labels to a graphics
object inserted using boxedeps.tex or BoxedArt.tex, make it a
firm rule to carefully adjust the bounding box using the trimming
commands of these packages, and also at least tentatively scale
and position the object. Beware of changing the grid inadvertently
after the labels have been positioned.  For example, correcting
the bounding box of a PostScript graphics object can foul up the
labels by changing the coordinate grid to which the labels are
attached. This is particularly true for the trimming  commands of
boxedeps.tex and BoxedArt.tex. However, as noted already, change
of scale is much less disruptive, and modest adjustments should be
well tolerated.

     Sometimes the labels protrude so far from the bounding box
of a figure that the figure has to be repositioned.  Best do this
by ad hoc spacing, say using \hglue and \vglue; altering the
bounding box would create a vicious circle.

     Remember that you are responsible for preventing labels
from overlapping. You are responsible for all label typography
including size and style. A label is really just about anything
that can be put in a TeX box. Note that spaces at the beginning
and end of labels will normally be suppressed; if you really want
them you must protect them with TeX braces.

     This package temporarily sets the \mathsurround parameter
of TeX to zero  while the labels are being affixed. This is done
because nonzero \mathsurround space would influence the position
of left and right aligned labels; then, when a texpert or printer
modifies mathsurround, diagram labeling might be disastrously
altered. There is a small price to pay involving labels that are
formatted as caption boxes including mathematics: you  may want or
need to specify an explicit mathsurround space within the caption
box; it will not influence anything outside.

     Those hostile to the use of * as separator between
the X and Y coordinates of label insertion points, are free to
impose another using \SetXYSeparator{<the new separator>}.  
Americans may prefer "," to "*" since they never use a 
comma as a decimal point; on the other hand, * may be more visible.

APPENDIX (I)  MERGING labelfig.tex LABELS INTO AN EPSF GRAPHICS OBJECT.

     As promised in the introduction, here is a recipe useful for
publishers. It works at least on Macintosh and at least for vectorized
graphics and Adobe type1 fonts.  (There is surely a similar recipe for
PCs under MSWindows.)

 (a)  Use boxedeps.tex utility to integrate the figure given by the eps
file, "x.eps" say, with a visible frame around it.  See
\ShowDisplacementBoxes command in boxedeps.tex.  To get precise results
automatically it is important to use the \Trim... commands of
boxedeps.tex making the "DisplacementBox" neatly fit the figure.

 (b)  Use the TeX printer driver and LaserWriter (versions >= 8.1.1) to
export to an EPSF the DVI page containing the integrated, labelled
figure. You now have an EPS file  "xx.eps"  that contains too much, and at
the wrong scale, and at wrong position.

 (c)  Convert the EPSF to an Adode Illustrator format EPSF using
the shareware utility called epsConvert by Sam Weiss
1993-- (currently $25).

 (d)  In Illustrator (or a compatible program), group the labels and the
"DisplacementBox"; copy them to the clipboard and paste them into "x.ps".
This step requires that all the label fonts be "visible to the Macintosh.

 (e)  Translate and scale the pasted group consisting of the labels plus
the "DisplacementBox" so as to make the "DisplacementBox" the bounding
box of (labelless) figure represented by "x.eps".  At this point the
labels will be correctly placed on the figure "x.eps".

 (f)  Ungroup and delete the "DisplacementBox".  The result is the
desired single EPS file, "x+.eps" say, It contains the original figure
plus its labels.  

     Using grouping and ungrouping appropriately in "x+.eps", a
publisher's staff can very efficiently improve label positions etc.

APPENDIX II)  SOME EXOTIC APPLICATIONS

     The grid of labelfig.tex is analogous to a light-table in
classical page makeup with wax or latex glue.  In principle, you
can use it to compose any page from its indivisible parts.  This
even has some of the artisanal charm of classical paste-up
provided you have a fast screen preview to make the process
"interactive".

     In practice labelfig.tex is a tool for nonstandard jobs.
Here are a few going beyond the labelling already discussed.

(I)  GRAPHICS INTEGRATION.

     This is accomplished by treating the imported graphics
objects as labels.  The underlying graphics object is then
typically an empty  \vbox to <dimension>{\vfill} in a TeX
\midinsert...\endinsert construction.  A label line
might be of the form

   (.1*.1) \\

The exact form of the special command varies from driver to
driver.  However, in the case of encapsulated PostScript graphics
(EPSF norm), by relying on boxedeps.tex, one can have the
following standard syntax (independant of driver  (see
boxedeps.doc for details.
  
  (.1*.1) \BoxedEPSF{MyFigure scaled <scale in mils>}\\

This may be slow since it requires TeX to read the PostScript
file to read bounding box using many complex macros.  So you
may want to try

  (.1*.1) \EPSFSpecial{MyFigure}{<scale in mils>}\\

which is fast and driver independant, but it squashes the
bounding box, normally to its lower left corner.

     Similarly for graphics of the Macintosh PICT norm ---
using BoxedArt.tex (same sources) in place of boxedeps.tex.

     This approach to integration is to be recommended when
one is assembling a composite graphics object.

 (II)  COMMUTATIVE DIAGRAM ENHANCEMENT

     Commutative diagrams or arrays of mathematical objects
connected by arrows of various sorts are common in mathematics.
The mathematical objects require the use of TeX.  Recently TeX
acquired a good collection of arrows of all slopes --- that of
LamSTeX --- plus pwerful macros to build the diagrams.

     However, even the LamSTeX collection is often
inadequate; it lacks for example double shafted arrows, dotted
arrows and curved arrows. Fortunately it is possible to produce
such arrows on an individual basis using sophisticated graphics
programs such as Illustrator and AldusFreehand (both serving
the EPSF norm) or using Metafont (with its public domain norm).
Since the creation of each new arrow is a work of love, you
probably want to limit the number of arrows by using LamSTeX
for most arrows. The 40K commutative diagram module of LamSTeX
has been adapted to work with AmSTeX and a copy may be posted
with LabelFig and related files. Unfortunately no one has yet
offered a version that works with Plain TeX or LaTeX.

       Suffice it here to say that when the exotic arrow has
been somehow imported into TeX, labelfig.tex treats it as a
label that one affixes to the commutative diagram.  Two other
steps will be treated in separate notes, namely the matter of
extracting the dimension specifications for the arrow and the
construction of the arrow --- for these steps are far from
unique and often depend intimately on your computer environment. 
Notes for the Macintosh-Textures-Illustrator combination are
found in the file ExoticArrows.doc.

 (III) NESTING 

Ingenuity pays off in exploiting labelfig.tex. One can
mix graphics and typography quite freely.  labelfig.tex is good
for freeform or overlapping arrangements, while boxedeps.tex (or
BoxedArt.tex) is best for regimented non-overlapping
arrangements --- and the two can be combined.

     The default behavior of labelfig.tex is not ideal 
for nesting objects, because to prevent trouble for beginners
the register for labels is globally cleared when \AffixLabels
concludes.  But there are switches available

      \LabelsGlobal      \LabelsLocal

which change this.  To understand this, extend the above test 
by something like:

 \LabelsLocal

 \SetLabels
    (.5*.5) AAA\\
 \endSetLabels

 {
 \SetLabels
    (.5*.5) ZZZ\\
 \endSetLabels
   \AffixLabels{\FirstQuadrant}
 }

   \AffixLabels{\FirstQuadrant}

     There are however potential pitfalls.  Neither
labelfig.tex nor boxedeps.tex has been tested under extreme
conditions. Problems may occur if their procedures are
indiscriminately nested. For boxedeps.tex (not labelfig.tex)
there is a precise cause for worry, namely many of its
variables are "global", which means that TeX braces will not
provide the protection one might expect.

COMMAND SUMMARY FOR labelfig.tex

  Here [...] means optional (one or zero)
       [...]* means any number of such constructs

  \SetLabels
    [[<P>](<X><Sep><Y>) <label> \\]*
  \endSetLabels
  \ShowGrid  
  \AffixLabels{<the figure>}

   --- <P> is tack position, one of eleven or empty
              order irrelevant

                   \L\T      \T      \R\T

                   \L\E      \E      \R\E

                     \L               \R

                   \L\B      \B      \R\B

   --- (<X><Sep><Y>) insertion point;
  <Sep> is separator, = * by default;
  \SetXYSeparator{<Sep>} changes it.
   <X> and <Y> are real numbers

  --- <label> a label to attach 

  --- <the figure> the figure to label 

  \GlobalLabels (default)     
  \LocalLabels  setting for nested constructs.

 \Grids makes ALL grids appear; \HideGrid then makes just next disappear.
 \noGrids returns to default.  The commands are always global.

 \GridLineWidth{<dimension>} adjusts width of grid lines. Default is very
small, to give "hairline" effect. If your grid lines are missing try
setting \GridLineWidth{1pt}.

 \Edges#1 globally changes the edge size of all grids to the numerical 
value #1, which must be 1, 10, 100, or 1000.  The default is 1.

VERSION HISTORY.
 --- Jan 1993: \Edges#1 and [??] option after \SetLabels
 --- July 1992: \Grids, \noGrids, \HideGrid;
       Gridlines become hairlines; \GridLineWidth{<dimension>}.
 --- Oct 1991, Jan 1992: \SetXYSeparator{<Sep>},  \LabelsGlobal,
       \LabelsLocal.
 --- July 1991: first release

Address for bugs and other feedback:

        Raymond S\'eroul
        IREM and Lab. de Typographie Informatise
        Univ. Rene Descartes
        Strasbourg

    Tel 33-88-41-63-45
    Email:  A18645@FRCCSC21.BITNET

        Laurent Siebenmann
        Mathematique, Bat. 425,
        Univ de Paris-Sud,
        91405-Orsay,
        France

    Tel 33-1-6941-7949; 
    Email: lcs@topo.math.u-psud.fr  

\newtheorem{theorem}{Theorem}[section]

\newtheorem{problem}[theorem]{Problem}
\newtheorem{conjecture}[theorem]{Conjecture}
\newtheorem{lemma}[theorem]{Lemma}
\newtheorem{proposition}[theorem]{Proposition}
\newtheorem{corollary}[theorem]{Corollary}
\theoremstyle{remark}
\theoremstyle{definition}
\newtheorem{definition}[theorem]{Definition}

\newcommand{\R}{\mathbb{R}}

\newcommand{\E}{\mathbb{E}}
\newcommand{\Z}{\mathbb{Z}}

\newcommand{\T}{\mathbb{T}}
\newcommand{\Py}{\mathsf{Py}}
\newcommand{\HH}{\mathbb{H}}
\newcommand{\F}{\mathbb{F}}

\newcommand{\fa}{\mathfrak{a}}
\newcommand{\fb}{\mathfrak{b}}

\newcommand{\calf}{\mathcal{F}}

\newcommand{\disco}{\mathsf{disco}}
\newcommand{\Good}{\mathsf{Good}}
\newcommand{\Prep}{\mathsf{Prep}}
\newcommand{\Avoid}{\mathsf{Avoid}}

\DeclareRobustCommand\longtwoheadrightarrow
     {\relbar\joinrel\twoheadrightarrow}

\numberwithin{equation}{section}
\numberwithin{figure}{section}

\let\qqed=\qed
\def\QED{\qqed\medskip}
\let\qed=\QED

\def \eps {\epsilon}

\def \P {{\bf P}}
\def\md{\mid}
\def\Bb#1#2{{\def\md{\bigm| }#1\bigl(#2\bigr)}}
\def\BB#1#2{{\def\md{\Bigm| }#1\Bigl(#2\Bigr)}}
\def\Bs#1#2{{\def\md{\mid}#1(#2)}}
\def\Pb{\Bb\P}
\def\Eb{\Bb\E}
\def\PB{\BB\P}

\def\Ps{\Bs\P}
\def\Es{\Bs\E}
\def\Pso#1{\Bs{\P_{#1}}}

\def \proof {{ \medbreak \noindent {\bf Proof.} }}
\def\proofof#1{{ \medbreak \noindent {\bf Proof of #1.} }}

\def\bl{\begin{lemma}}
\def\el{\end{lemma}}
\def\bth{\begin{theorem}}
\def\eth{\end{theorem}}
\def\bc{\begin{corollary}}
\def\ec{\end{corollary}}
\def\bcj{\begin{conjecture}}
\def\ecj{\end{conjecture}}
\def\bpr{\begin{proposition}}
\def\epr{\end{proposition}}
\def\bde{\begin{definition}}
\def\ede{\end{definition}}
\newcommand{\be}{\begin{eqnarray}}
\newcommand{\ee}{\end{eqnarray}}
\newcommand{\bes}{\begin{eqnarray*}}
\newcommand{\ees}{\end{eqnarray*}}

\usepackage{mathrsfs}

\def\1{1\!\! 1}
\def\cL{\mathcal L}
\def\cT{\mathcal T}
\def\cA{\mathcal A}
\def\cB{\mathcal B}
\def\cF{\mathcal F}
\def\cC{\mathcal C}
\def\cH{\mathcal H}
\def\cD{\mathcal D}
\def\cG{\mathcal G}

\def\calt{\mathcal T}

\def\cM{\mathcal M}

\def\UST{\mathsf{UST}}
\def\FUSF{\mathsf{FUSF}}
\def\WUSF{\mathsf{WUSF}}
\def\LERW{\mathsf{LERW}}
\def\FMSF{\mathsf{FMSF}}

\title{The Free Uniform Spanning Forest is disconnected\\ 
in some virtually free groups, depending on the generator set}
\author{
G\'abor Pete \and \'Ad\'am Tim\'ar
}
\date{January 24, 2021}

\begin{document}
\maketitle

\begin{abstract} 
We prove the rather counterintuitive result that there exist finite transitive graphs $H$ and integers $k$ such that the Free Uniform Spanning Forest in the direct product of the $k$-regular tree and $H$ has infinitely many trees almost surely. 

This shows that the number of trees in the $\FUSF$ is not a quasi-isometry invariant. Moreover, we give two different Cayley graphs of the same virtually free group such that the $\FUSF$ has infinitely many trees in one, but is connected in the other, answering a question of Lyons and Peres \cite{LPbook} in the negative.


A version of our argument gives an example of a non-unimodular transitive graph where $\WUSF\not=\FUSF$, but some of the $\FUSF$ trees are light with respect to Haar measure. This disproves a conjecture of Tang \cite{Pengfei}.
\end{abstract}

\section{Intro}

The Free Uniform Spanning Forest $\FUSF$ is one of the most standard random spanning forests of infinite graphs, obtained as the weak limit of the uniform random spanning trees $\UST$ in any exhaustion of the infinite graph by finite subgraphs. In any transitive graph, its law is invariant under the automorphisms of the graph. It may be regarded as the Free $\mathsf{FK}(p,q)$ random cluster model with $q=0$ at its critical point $p=0$, it is a determinantal process, and is especially interesting due to its connections to measurable group theory: in any Cayley graph of a group $\Gamma$, its expected degree is $2+2\beta^{(2)}_1(\Gamma)$, where $\beta^{(2)}_1(\Gamma)$ is the {\it first $\ell^2$-Betti number} of the group, the von Neumann dimension of the space of harmonic functions of finite Dirichlet energy. In particular, we have the equality $\FUSF=\WUSF$ with the Wired Uniform Spanning Forest if{f} $\beta^{(2)}_1(\Gamma)=0$. See \cite{BLPS:USF} and \cite[Chapter 10]{LPbook} for thorough studies of the  $\FUSF$; some more recent papers are \cite{HuNa,indist,AHNR,HuNaPlanar}.

We will mostly work in the direct product graph $\T^k \times H$, where $\T^k$ is the $k$-regular infinite tree with $k\ge 3$, while $H$ is a finite vertex-transitive graph. 
Typical examples are the product Cayley graphs of the virtually free groups $\F_r \times \Gamma$, where $\F_r$ is a free group on $r\ge 2$ generators and $\Gamma$ is a finite group. The $\FUSF$ on some tree-like graphs was recently studied, among other topics, in \cite{Pengfei}. In particular, Tang proved that, for any $k$, the $\FUSF$ in $\T^k\times \Z_2$ (where $\Z_2$ is the path on 2 vertices, i.e., a single edge) is connected almost surely; this was later generalized in \cite{ABIT} for an arbitrary fixed weight on the $H$-edges. Tang made the innocent-looking conjecture that the connectedness holds more generally, for the direct product $\T^k \times H$ with any $k\ge 3$ and any finite transitive graph $H$ (no edge weights).  See Remark 5.9 in that paper. 
 Here we are disproving this conjecture.

\begin{theorem}[Disconnected $\FUSF$]\label{t.disco}
For every $d$ there is $k_d$ such that if $\T^k$ is the $k$-regular infinite tree with $k\geq k_d$, and $H$ is a connected  finite $d$-regular transitive graph on more than $k^{5/2}$ vertices, then the $\FUSF$ of $\T^k \times H$ is disconnected almost surely. In fact, it has infinitely many components.
\end{theorem}

One striking corollary of our result is that the number of trees in the $\FUSF$ of Cayley graphs is {\it not a quasi-isometry invariant}, as opposed to several similar properties: the number of trees in the $\WUSF$ \cite[Corollary 10.25]{LPbook}, the property $\WUSF\not=\FUSF$ \cite{Soardi, BLPS:USF}, or equivalently, the infinite-endedness of all the $\FUSF$ trees (the equivalence follows from \cite{recurr} and \cite{HuNa,indist}). (Note, nevertheless, that without transitivity of the graph the number of components is not a quasi-isometry invariant even when $\WUSF=\FUSF$, as the example in \cite{Itai} shows.) 
With some extra work, we prove here that the number of trees is not even the same for different Cayley graphs of a fixed group (even though the expected degree of the $\FUSF$ depends only on the group, because of the connection to $\beta^{(2)}_1(\Gamma)$). This answers a question of  Lyons and Peres \cite[Question 10.50]{LPbook} in the negative:

\begin{theorem}[Dependence on the generating set]\label{t.genset}
For $k$ large enough, the group $\F_k \times \Z_{k^9}$ (the direct product of a free group and a cyclic group) has a Cayley graph (the direct product of the tree $\T^{2k}$ and the cycle $C_{k^9}$) in which the $\FUSF$ has infinitely many components, and another Cayley graph (the direct product of the tree $\T^{2k}$ and the complete graph $K_{k^9}$) in which the $\FUSF$ is connected. 
\end{theorem}

To our knowledge, this is the first instance of a standard statistical physics model that shows such {\it non-universal critical behavior}.

Another corollary of Theorem~\ref{t.disco} is that although the $\FUSF$ might be connected in every quasi-transitive (or more generally, unimodular random) planar graph (see \cite{AHNR} for a large subclass), this for sure cannot be extended from planar graphs to an arbitrary minor-closed family. This also means that a positive answer to \cite[Question 8]{soficplanar}, extending treeability and soficity of unimodular random graphs from the planar case to graphs with arbitrary excluded minors, cannot be done via the strategy of \cite{AHNR}, using the $\FUSF$.

It should be mentioned that \cite[Question 11.37]{LPbook} asks whether the Free {\it Minimal} Spanning Forest $\FMSF$ is connected in any graph that is roughly isometric to a tree. A key difference from our situation is that \cite[Theorem 1.3]{MSF} says that the union of $\FMSF$ with an independent $\mathsf{Bernoulli}(\eps)$ bond percolation is always connected, for any $\eps>0$. This is something that we do not know for the $\FUSF$ in our graphs, which also brings us to our next remark.

A well-known question of Damien Gaboriau \cite{Gaboriau} is whether the so-called {\it measurable cost} of any group $\Gamma$ is equal to $1+\beta^{(2)}_1(\Gamma)$. He pointed out (see \cite[Question 10.12]{LPbook}) that a positive answer would follow if, in every Cayley graph and any $\eps>0$ there was a connected invariant bond percolation $\omega$ that contains $\FUSF$, but $\omega\setminus\FUSF$ has density at most $\eps$. Interesting examples are the infinite Kazhdan groups: here $\beta^{(2)}_1(\Gamma)=0$, hence $\WUSF=\FUSF$, by \cite{BekVal}; thus non-amenability together with  \cite[Theorem 13.7]{BLPS:USF} imply that adding an independent $\mathsf{Bernoulli}(\eps)$ bond percolation does not work; on the other hand, adding some much trickier invariant percolation does work \cite[Remark 2.2]{costKazhdan}. In the examples of our Theorem~\ref{t.disco}, we have $\WUSF\not=\FUSF$ (because transitive graphs with infinitely many ends have harmonic functions with finite Dirichlet energy), so it is tempting to speculate that they could provide a negative answer to Gaboriau's question. However, we have been unable to prove anything in this direction. In particular, it remains open if any two trees in our $\FUSF$ touch each other at finitely many places, similarly to Bernoulli percolation \cite{nontouch} or $\WUSF$ clusters in $\Z^d$ with $d\ge 9$ \cite{4812}.

Let us note that for any {\it infinite} transitive graph $\HH$, it has been known for long \cite{BLPS:USF} that the $\FUSF$ of $\T^k \times \HH$ has infinitely many components. More generally, in the direct product of any non-amenable transitive graph with any infinite transitive graph, there is no invariant probability measure on the set of subtrees \cite{PePe} (even without the requirement of being spanning trees). However, for any {\it finite} transitive graph $H$, a uniform random translate of $\T^k$ gives an invariant random subtree, hence a general non-treeability argument could not imply our Theorem~\ref{t.disco}. In fact, all disconnectedness results on the $\FUSF$ that we know of have been obtained so far either by proving that $\WUSF=\FUSF$ and knowing that the $\WUSF$ trees are small (e.g., recurrent, hence one or two-ended \cite{recurr}); or by noticing that even when $\WUSF \not=\FUSF$, the $\FUSF$ may be similar to the $\WUSF$, as in the free product $\Z^5 * \Z_2$; or by a general non-treeability result, which applies not only to the $\FUSF$ but to any invariant spanning forest. In contrast, our proof is in a treeable group, specific to the $\FUSF$, in a situation where $\WUSF\not=\FUSF$. The reason for having no earlier $\FUSF$-specific results is that this is quite a mysterious object: while the $\WUSF$ can be generated in infinite graphs directly by Wilson's algorithm rooted at infinity, using loop-erased random walks \cite{Wilson}, or by the Interlacement Aldous-Broder algorithm \cite{IAB}, no such method is known for the $\FUSF$. Indeed, we will use Wilson's algorithm in finite balls of the graph, then take the limit.

As a follow-up to the present paper, the preprint \cite{ABIT} studies the $\FUSF$ on direct products $\T^k \times H$ with edge weights $c \in (0,\infty)$ for the edges of $H$. It is proved there that for any $\T^k \times H$, if $c$ is large enough, then the $\FUSF$ is connected. This would immediately imply the connectedness direction of our Theorem~\ref{t.genset} if we allowed for weighted generating sets. However, getting a standard unweighted generating set has some value: e.g., for Bernoulli percolation on nonamenable groups, it was proved by \cite{PSN} that some weighted generating set has $p_c<p_u$, and it took fifteen years to achieve the same result without weights \cite{Thom}. 
\medskip

A version of our construction gives a counterexample to Conjecture 1.2 of \cite{Pengfei}, in a strong way. A transitive graph $G$, with full automorphism group $\Gamma$, is called {\it unimodular} if, for every pair of neighbors $x,y$, we have $|\Gamma_x y| = |\Gamma_y x|$, where $\Gamma_x=\{\gamma\in\Gamma: \gamma(x)=x\}$ is the stabilizer subgroup, and $\Gamma_x y=\{\gamma(y) : \gamma\in \Gamma_x\}$ is the orbit of $y$. For instance, every Cayley graph is unimodular. See \cite[Chapter 8]{LPbook} on background on unimodularity and its connections to invariant percolations. For non-unimodular transitive graphs, it is worth looking at an invariant Haar-measure $\mu$ on the locally compact $\Gamma$, which gives finite but non-equal weights to the stabilizers: 
$$
\frac{\mu(\Gamma_x)}{\mu(\Gamma_y)} = \frac{|\Gamma_x y|}{|\Gamma_y x|},
$$ 
for any $x,y \in V(G)$. A subset $\cC \subset V(G)$ is called {\it light} if $\sum_{x\in\cC} \mu(\Gamma_x) < \infty$. It was proved in \cite[Theorem 1.1]{Pengfei} that the trees of $\WUSF$ in any non-unimodular transitive graph (and more generally, whenever there is a closed non-unimodular subgroup of automorphisms that acts transitively on $G$) are light. His Conjecture 1.2 stated that the opposite holds for $\FUSF$, when $\WUSF\not=\FUSF$. Since our examples in Theorem~\ref{t.disco} do have transitive closed non-unimodular subgroups (the automorphisms fixing an end of the tree), they already give counterexamples to the more general conjecture. Nevertheless, with a bit of more work, we can also give counterexamples where the full automorphism group is non-unimodular. Note here that there is a usual way of producing a non-unimodular transitive graph from a graph with a non-unimodular transitive subgroup of automorphisms by adding some edges in a transitive way, as in the grandmother graph; however, since we have already seen that the number of $\FUSF$ components is not a quasi-isometry invariant, it is unclear what the effect of such a ``small'' change would be.

\begin{theorem}[Non-unimodular lightness]\label{t.light}
There exists a non-unimodular transitive graph $G$ in which $\WUSF\not=\FUSF$, but $\FUSF$ has some light clusters. 
\end{theorem}

The second part of \cite[Conjecture 1.2]{Pengfei} was that, for nonunimodular transitive graphs, $\WUSF\not=\FUSF$ implies that all the trees of $\FUSF$ have branching number larger than 1. (True in the unimodular case, because the average degree being strictly larger than 2 implies invariant non-amenability \cite[Section 8]{urn}.) Our construction is not a counterexample to this conjecture.
\medskip

The dis/connectedness results discussed above give rise to a nontrivial graph parameter: for any finite graph $H$ we let
$$
\disco(H):=\min\big\{k : \FUSF(\T^k \times H) \textrm{ is disconnected}\big\} \in \{3,4,\dots,\infty\}.
$$
The earlier results on the connectedness of $\FUSF$ in $\T^k\times P_2$ say that $\disco(P_2)=\infty$. Our Theorem~\ref{t.disco} implies that if $\ell$ is large enough, then the cycle $C_\ell$ of length $\ell$ has $\disco(C_\ell) < \infty$. Several specific open questions on this graph parameter are discussed in~Section~\ref{s.open}.
\medskip

To conclude this introduction, let us say a few words about our proof strategies and the organization of the paper.

The ball of radius $n$ around a fixed root $o\in\T^k$ will be denoted by $T_n$, while the sphere of radius $n$ will be denoted by $S_n$. We will generate the $\UST$ in $T_n \times H$  by Wilson's algorithm \cite{Wilson,LPbook}, first taking the loop-erased random walk $\LERW$ from ${\fa}=(o,h_{\fa})$ to ${\fb}=(o,h_{\fb})$, where $h_{\fa}\not= h_{\fb} \in H$ are arbitrary. See Section~\ref{s.product} for the definitions. In the setting of Theorem~\ref{t.disco}, we will prove that the $\LERW$ from $\fa$ to $\fb$, with a positive probability that does not depend on the radius $n$, will contain some close-to-the-boundary vertex $(\tilde z,h_{\tilde z}) \in S_{n-7} \times H$. This will easily imply the theorem. Finding such a  $(\tilde z,h_{\tilde z})$ will go as follows. 

We will find that the simple random walk trajectory from ${\fa}$ to ${\fb}$ with uniformly positive probability hits a ``bag" $\{z\}\times H$ with $z \in S_n$ in such a way that the part of the trajectory before hitting $\{z\}\times H$, denoted by $\pi_\mathrm{there}$, and the second part $\pi_\mathrm{back}$ after leaving $\{z\}\times H$ intersect each other only outside $T_{n-8} \times H$. Then $\tilde z$ will be the ancestor of $z$ in $S_{n-7}$. To find such a $z$, we will have to make sure that there are no intersections in either of the following ways: {\bf (1)} outside the ray of bags between $\{o\}\times H $ and $\{\tilde z\}\times H $; {\bf (2)} in some bag of this ray. 

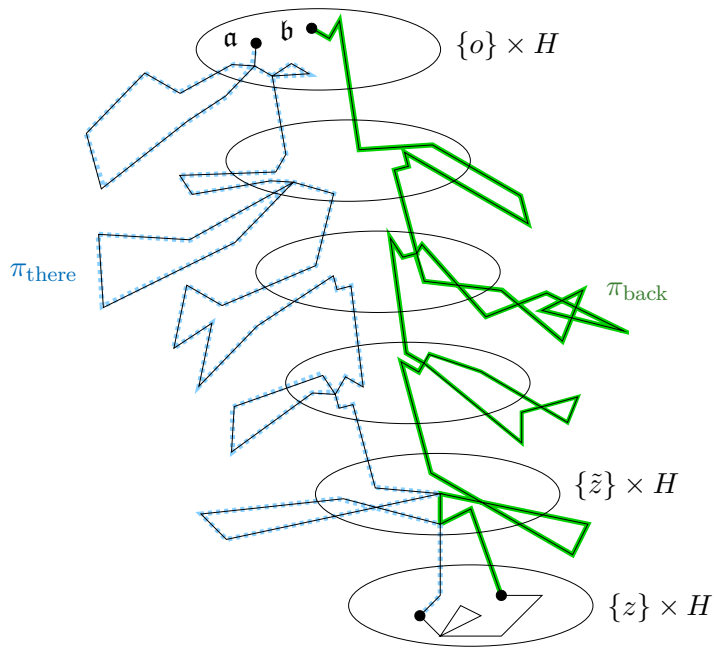
\begin{figure}[htbp]
\SetLabels
(0.28*0.94) $\fa$\\
(0.38*0.95) $\fb$\\
(0.77*0.93)  $\{o\}\times H$\\
(0.98*0.24) $\{\tilde{z}\}\times H$\\
(1.04*0.05) $\{z\}\times H$\\
(-0.05*0.58) \textcolor{NavyBlue}{$\pi_\mathrm{there}$}\\
(1*0.55) \textcolor{OliveGreen}{$\pi_\mathrm{back}$}\\
\endSetLabels
\centerline{
\AffixLabels{
\includegraphics[width=0.44\textwidth]{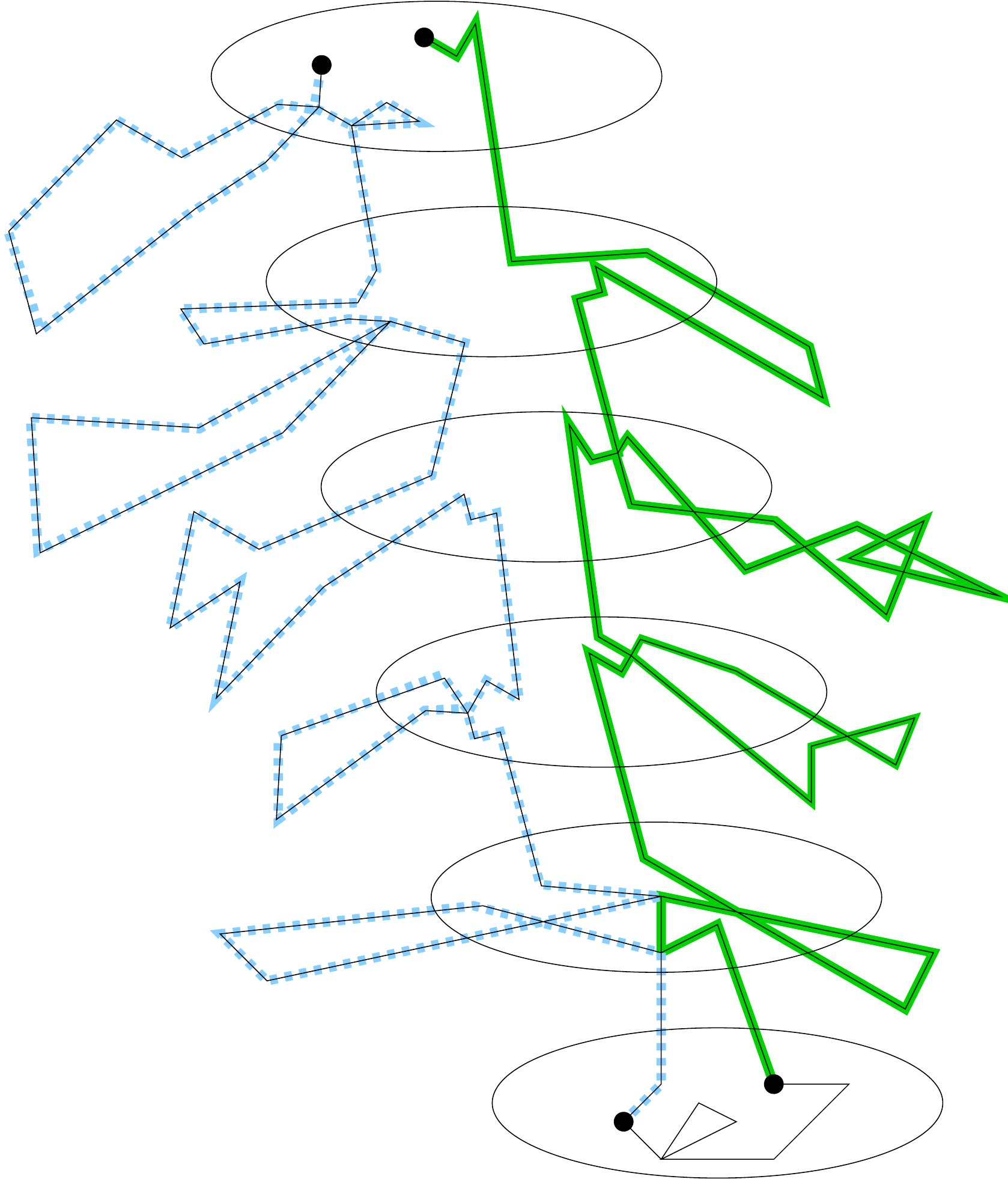}
}}
\caption{Strategy of proof in $T_n \times H$. In this picture, $\pi_\mathrm{there}$ and $\pi_\mathrm{back}$ do not intersect each other inside $T_{n-2}$, hence we obtain $\tilde z \in S_{n-1}$ such that the $\LERW$ path from ${\fa}$ to ${\fb}$ intersects $\{\tilde z\} \times H$.}
\label{f.strategy}
\end{figure}


To guarantee (1), and also to help with (2), we will ensure that neither $\pi_\mathrm{there}$ nor $\pi_\mathrm{back}$ makes any backtracking on the ray of bags from $o$ to $z$, and furthermore, the bags that $\pi_\mathrm{there}$ enters outside this ray are different from the bags that $\pi_\mathrm{back}$ enters. See Figure~\ref{f.strategy}. These requirements concern only the tree-coordinate of the random walk, and it is indeed possible to find $z$ such that they (and hence (1)) are satisfied, as we will prove in Proposition~\ref{p.localtime}.

It remains to rule out intersections as in (2). Here the $H$-coordinates will play the main role. The intuition is that the visits in a typical bag $\{v\}\times H$ are not too long (since $k$ is large compared to $d$), and the places where the walker enters $\{v\}\times H$ from the outside are likely to be far from each other, because these entrances tend to be separated by long time intervals (until the walk on the tree returns) and because $H$ is large. To elaborate this argument will require some work, presented in Section~\ref{s.product}.
\medskip

For Theorem~\ref{t.genset}, the idea is to start with a small degree $d$ but large $H$ compared to $k$, so that Theorem~\ref{t.disco} applies, then change the generating set so that we get the complete graph on $H$. This makes the random walk that generates the $\LERW$ spend a lot of time in each bag $\{v\}\times H$ before  moving in the tree-coordinate, making it very likely that the loop-erasure erases every long excursion away from the root bag $\{o\}\times H$. The details are worked out in Section~\ref{s.change}.
\medskip

In Section~\ref{s.nonunimod}, we prove Theorem~\ref{t.light} on lightness in the non-unimodular setting. Here the task is to modify the tree-proof to a well-chosen non-unimodular transitive graph, then argue that there are infinitely many components in the $\FUSF$, which makes at least some of them light.
\medskip

We conclude the paper with several open problems in Section~\ref{s.open}, including the ones on Gaboriau's question and on our new graph parameter $\disco(H)$ for finite graphs $H$.


\section{Born to be alive}\label{s.tree}


A key observation about the tree-coordinate of the random walk will be the following proposition, somewhat interesting in its own right. Consider simple random walk $(Y_t)_{t\geq 0}$ on $T_n$, started at the root $o$, until the first return time $\tau_o^+:=\min\{t > 0 : Y_t =o\}$. 

\bpr[Viable rays]\label{p.localtime}
For any $k$ large enough, with a positive probability that may depend only on $k$, there is a $z\in S_n$ in $\T^k$ such that, denoting the ray from $o$ to $z$ in $T_n$ by $\gamma=(o=\gamma_0,\gamma_1,\dots,\gamma_n=z)$, we have:
\begin{itemize}
\item all the edges on the ray $\gamma$  are crossed exactly twice until $\tau_o^+$ (once on the way from $o$ to $z$, once on the way back);
\item on the way from $o$ to $z$, for every $i=1,\dots,n-1$, the number of excursions away from $\gamma_i$ before taking the edge  $(\gamma_i,\gamma_{i+1})$ is at most $k/2$;
\item denoting by $E_i$ and $F_i$ the set of edges incident to a vertex $\gamma_i$ but not on $\gamma$ that are crossed on the way to $z$, and on the way back from $z$ to $o$, respectively, we have that $E_i\cap F_i = \emptyset$ for all $i=1,\dots,n-1$. 
\end{itemize}
\epr

Such a ray typically has the property that all its vertices have positive but small local times (of order $k$) until $\tau_0^+$. It is possible that, using the Dynkin isomorphism theorem \cite{Dynkin}, such a result could be proved via the Gaussian Free Field on $\T^k$; see \cite{DLP, Lupu, Zhai} for such arguments. However, since we also need the more refined statement on the edges incident to the ray, we have not tried to make this connection precise. Let us emphasize that a typical ray to $S_n$, or the first ray along which we reach $S_n$, do not satisfy the proposition; we have to work to find such rays.

\proofof{Proposition~\ref{p.localtime}} Pick a leaf $z\in S_n$, denote the ray from $o$ to $z$ by $\gamma=\gamma(z)$, the stopping times $\tau_z:=\min \{ t : Y_t=z\}$ and $\tau^+_o$ as before, and define the events
\be\label{e.AL}
\begin{aligned}
\cA_z &:= \Big\{ \textrm{the edge }(\gamma_{i-1},\gamma_{i})\textrm{ is crossed exactly twice by }(Y_t)_{t=0}^{\tau^+_o} \textrm{, for all } i=1,2,\dots,n\Big\},\\
\cL_z &:= \Big\{ \big|\big\{t\in \{1,\dots,\tau_z\} : Y_t=\gamma_i \big\} \big| \leq k/2+1 \textrm{, for all } i=1,2,\dots,n\Big\}.
\end{aligned}
\ee
Furthermore, let $E_i$ and $F_i$ be the set of edges as defined in Proposition~\ref{p.localtime}, and 
define the event
\be\label{e.avoid}
\cB_z := \cA_z \cap \cL_z \cap \big\{E_i \cap F_i = \emptyset \textrm{ for all } i=1,2,\dots,n-1 \big\}.
\ee

Let us calculate $\Ps{\cB_z}$. The first step has to be $\Ps{Y_1 = \gamma_1}=1/k$, and then, for each $\gamma_i$, $ i=1,\dots,n-1$, the walk $(Y_t)$ may take excursions away from $\gamma_i$, but it has to choose $\gamma_{i+1}$ before $\gamma_{i-1}$, and the number of excursions has to be at most $k/2$. The probability of this event $\mathsf{There}_i$, with the extra condition that there are precisely $j \ge 0$ excursions, is 
\be\label{e.There}
\Pb{\mathsf{There}_i, \textrm{ with }j\textrm{ excursions}} = \left(1-\frac2k\right)^j \frac2k \, \frac12, 
\ee
independently of what happens at other $\gamma_i$'s. When we arrive at $Y_{\tau_z}=z$, we have already sampled the edge sets $E_i$, $i=1,\dots,n-1$. Then, at each $\gamma_i$, for $i=n-1,n-2,\dots, 1$, we have to choose  $\gamma_{i-1}$ before $\gamma_{i+1}$, an event we will denote by $\mathsf{Back}_i$; furthermore, the excursions away from $\gamma_i$ have to produce an edge set $F_i$ that is disjoint from $E_i$. The probability of everything together, independently of $i$, is 
\be\label{e.pk}
p_k:= \Pb{\mathsf{There}_i,\ \mathsf{Back}_i, \textrm{ and } E_i \cap F_i = \emptyset } 
\ge \sum_{j=0}^{\lfloor k/2\rfloor} \left(1-\frac2k\right)^j\frac1k \frac{1}{j+2} 
\asymp \frac{\log k}{k}\,,
\ee
because if we have $j\leq k/2$ excursions in~(\ref{e.There}), then $|E_i| \leq j$, thus the walk on the way back has to avoid at most $j+1$ neighbors before choosing $\gamma_{i-1}$ (the edges of $E_i$ plus the edge to $\gamma_{i+1}$), which has success probability $1$ in at most $j+2$ . The asymptotic formula at the end simply follows from the exponential factor being between 1 and $1/e$ for all $0 \leq j \leq k/2$; the symbol $\asymp$ means ``up to positive universal constant factors'', independently of $k$ or $n$. 

The events of (\ref{e.pk}) for different $i$'s are independent from each other, hence we have
\be\label{e.Bz}
\Ps{\cB_z} = \frac{1}{k} \, p_k^{n-1}.
\ee

Let $Z_n$ be the set of leafs $z\in S_n$ that satisfy the event $\cB_z$. Then we have the first moment
\be\label{e.Z1}
\E |Z_n| = k (k-1)^{n-1} \, \frac{1}{k} \, p_k^{n-1},
\ee
which goes to infinity as $n\to\infty$ if $k$ is large enough, by~(\ref{e.pk}).

To estimate the second moment $\E |Z_n|^2$, let $z,v \in S_n$ be leafs such that their last common ancestor is $w\in S_m$, with $m\ge 1$. We claim that
\be\label{e.BzBv}
\Pb{\cB_z \cap \cB_v} \asymp_k p_k^{2n-m} \,,
\ee
where $p_k$ is defined in (\ref{e.pk}), and $\asymp_k$ means ``up to constant factors that may depend on $k$, but not on $n$ or $m$''. 

Indeed, the first step in $(Y_t)$ has to be towards $w$; then we have to reach $w$ without ever stepping backwards along the ray from $o$ to $w$; then we have to step towards $z$ or $v$ before stepping backwards towards $o$; then we have to reach the chosen leaf without backward moves; then we have to go back to $w$ without backward moves, and with the $F_i$ sets avoiding the $E_i$'s; at $w$, we have to step towards the other leaf before stepping towards $o$; if we define $F'_m$ to be the set of edges emanating from $w$ that are crossed after reaching $w$ after the first leaf, but before the step towards the second leaf, we must have $E_m\cap F'_m=\emptyset$; from $w$ we have to reach the other leaf without ever stepping backwards; then we have to go back to $w$, without backward moves, and with the $F_i$ sets avoiding the $E_i$'s also along this branch; at $w$, we have to move towards $o$ before moving towards $z$ or $v$ again, and the edge set $F''_m$ produced by the excursions before that has to be disjoint both from $E_m$ and $F'_m$; then we have to reach $o$ without ever stepping backwards, again with the $F_i$ sets avoiding the $E_i$'s. We have $m-1+2(n-1-m)=2n-m-3$ of these conditions at vertices other than $w$, each with success probability $p_k$, independently from each other. At $w$, the conditions are obviously possible to satisfy  if $k\ge 3$ (at the first visit go straight towards $z$, at the second visit go straight towards $w$, at the third visit go straight towards $o$), happening with probability at least $1/k^3$ and at most 1. So, the probability altogether is between $p_k^{2n-m-3} / k^3$ and $p_k^{2n-m-3}$, which can be written as~(\ref{e.BzBv}).

By going through all possible last common ancestors $w$, from~(\ref{e.BzBv}) we get
\be\label{e.Z2}
\E |Z_n|^2 = \sum_{z,v \in S_n} \Pb{\cB_z \cap \cB_v}
\asymp_k \sum_{m=1}^n k(k-1)^{m-1}\, (k-1)^{2(n-m)}\, p_k^{2n-m} \asymp_k \big((k-1)p_k\big)^{2n}\,,
\ee
if $k$ is large enough, since $(k-1) p_k \to \infty$ holds by~(\ref{e.pk}), hence the $m=1$ term will dominate.

Comparing (\ref{e.Z1}) and (\ref{e.Z2}), the Cauchy-Schwarz second moment method gives us 
$$
\Ps{|Z_n|>0} \geq \frac{(\E |Z_n|)^2}{\Es{ |Z_n|^2}} \asymp_k 1,
$$
finishing the proof of Proposition~\ref{p.localtime}.
\qed

\section{Stayin' alive}\label{s.product}

In this section, we will first recall how to generate the $\FUSF$ via an exhaustion by finite graphs and the loop-erased random walk $\LERW$ inside each finite graph. Then we will consider the random walk in $\T^k \times H$, together with its projection to $\T^k$, and show that some of the viable rays found in Section~\ref{s.tree} correspond to trajectories in the product graph that survive the loop-erasure, provided that $H$ is large enough (compared to $k$). This way, we get distinct paths in the $\FUSF$ from two neighboring vertices to infinity.

As we briefly explained in the Introduction, the $\FUSF$ of an infinite graph $G$ is defined as the weak limit of the sequence $\UST(G_n)$, where $(G_n)_{n\ge 1}$ is any increasing sequence of connected finite subgraphs of $G$ such that $\bigcup_{n\ge 1} G_n = G$, and $\UST$ is the uniform measure on all spanning trees of the finite graph. The limit exists and is independent of the sequence $(G_n)$ by some electric network monotonicity arguments \cite[Chapter 10]{LPbook}. On a connected finite graph $G$, we can use the {\it loop-erased random walk} $\LERW$ to construct $\UST(G)$ with Wilson's algorithm \cite{Wilson}. Choose two vertices $x_0, x_1$ of $G$, and produce a simple path from $x_1$ to $x_0$ by running a random walk from $x_1$ until hitting $x_0$, and erasing all cycles created by the trajectory, in the order of creation. Then pick some $x_2$, start a walk from $x_2$ until we hit the path between $x_0$ and $x_1$, take the loop-erasure of it, and so on, always walking from $x_i$ until we hit the already existing tree, repeating until all the vertices become part of the tree.

Our infinite graph will be a direct product $G=\T^k\times H$, often denoted by $\T^k\square H$,  where the vertex set is just the set of pairs, and the neighbors of $(t,h)$ are the vertices $(t',h)$ with $\{t,t'\}\in E(\T^k)$ and the vertices $(t,h')$ with $\{h,h'\}\in E(H)$.

\proofof{Theorem~\ref{t.disco}} We take the exhaustion $G_n = T_n \times H$ of $G=\T^k\times H$. Our first step in Wilson's algorithm is to take the loop-erased random walk $\LERW$ from ${\fa}=(o,h_{\fa})$ to ${\fb}=(o,h_{\fb})$, where $h_{\fa}\not= h_{\fb} \in H$ are arbitrary. We will prove that the $\LERW$ from $\fa$ to $\fb$, with a  probability greater than a positive number $p$ that does not depend on $n$, will contain some close-to-the-boundary vertex $(\tilde z,h_{\tilde z}) \in S_{n-7} \times H$. Then, for any fixed finite subgraph $U$ of $G$, if $n$ is large enough so that $G_n$ contains $U$, but $S_{n-7}\times H$ is already disjoint from $U$,  and the above event for the $\UST(G_n)$-path between $\fa$ and $\fb$ occurs, then the intersection of this $\UST(G_n)$-path with $U$ will not connect $\fa$ and $\fb$. Hence, in the weak limit as $n\to\infty$, the $\FUSF$-component of $\fa$ will be different from the component of $\fb$ with probability at least $p$. 

One way to complete the proof from here is that the number of trees in the $\FUSF$ in any unimodular transitive graph was shown in  \cite{indist} and \cite{HuNa} to be either one a.s., or infinite a.s., hence we have to be now in the second case. We will also give a direct proof for our very special product graph, immediately extendable to the non-unimodular graph of Section~\ref{s.nonunimod}, via Wilson's algorithm, at the end of this section. And, we will give yet another proof, using the Mass Transport Principle, which again works both for unimodular and non-unimodular transitive graphs that are tree-like in some sense, and in a larger generality than the $\FUSF$, in Proposition~\ref{p.inftymany}. Some readers might prefer the more specific Wilson's algorithm proof, some readers might prefer the more general MTP proof, but in any case, not relying on unimodularity will be important for Theorem~\ref{t.light}.

We now turn to the study of the $\LERW$ from $\fa$ to $\fb$. The random walk on $G_n$ from $\fa$ to $\fb$ for which we apply the loop-erasure will be denoted by $(X_t)_{t\geq 0}$. The first coordinate of $(X_t)_{t\geq 0}$ is a lazy random walk on $T_n$, denoted by $(Y_t)_{t\geq 0}$. As before, we fix $z\in S_n$, and let $\tau_z$ and $\tau^+_o$ denote the hitting times for the projection $(Y_t)$. Condition on  the event $\cB_z$ of~(\ref{e.avoid}), but with the local times at $\gamma_i$ in the definition of $\cL_z$ being understood as the number of ``essentially different visits'', i.e., with the lazy steps removed from $(Y_t)$. The last time before $\tau_z$ that $(X_t)$ is in $\gamma_i \times H$ is denoted by $\alpha_i$, and the first time after $\tau_z$ that $(X_t)$ is in $\gamma_i \times H$ is denoted by $\beta_i$. (For $i=n$, we mean $\alpha_n = \beta_n = \tau_z$.) Furthermore, the number of actual (non-lazy) steps until $\alpha_i$ in the $\T^k$ and $H$ coordinates will be denoted by $\alpha^\T_i$, $\alpha^H_i$, respectively, and similarly for $\beta_i$. The first ingredient in our proof will be that, with a uniformly positive probability, quite a long time passes between entering consecutive bags of $\gamma\times H$ during $\pi_\mathrm{back}$.

\begin{lemma}\label{l.back}
Let $\cF(\beta_i)$ be the sigma-algebra generated by $(X_t)_{t=0}^{\beta_i}$. Then, for any $1\leq i \leq n-7$, 
$$
\Pb{\beta^H_{i-1}-\beta^H_i > k^5 \md \cB_z, \cF(\beta_i)} >  b
$$ 
for any large enough $k$, with a constant $b>0$ that does not depend on $i$, $n$, or $k$.
\end{lemma}

\proof Since we are conditioning on an event concerning the entire random walk trajectory, $\cB_z$, we have to be careful what the exact effect of this is. Namely, for any such event $\cB$, the original random walk transition probabilities get reweighted by a {\it Bayesian factor}:
\be\label{e.Bayes}
\PB{X_{t+1} \md (X_s)_{1\leq s \leq t}, \cB } = \Pb{ X_{t+1} \md  (X_s)_{1\leq s \leq t} } \frac{\Pb{\cB \md  (X_s)_{1\leq s \leq t+1}}}{ \Pb{\cB \md  (X_s)_{1\leq s \leq t}}}.
\ee

Now, for the lemma, it is enough to prove that 
\be\label{e.backT}
\Pb{\beta^\T_{i-1}-\beta^\T_i > 2 k^6 \md  \cB_z, \cF(\beta_i) } >  b',
\ee 
with some constant $b'>0$, by the following reasoning. Whenever $X_t \in \gamma_i \times H$ at some time $t\ge \beta_i$, the conditioning on $\cB_z$ forbids the $\T^k$-steps through $(\gamma_i,\gamma_{i+1})$ and $E_i$, while the other $\T^k$-steps and all the $H$-steps are available. In other words, the Bayesian factor from~(\ref{e.Bayes}), with $\cB=\cB_z$ and also conditioned on $\cF(\beta_i)$, is zero for $X_{t+1} \in \gamma_{i+1} \times H$ and for $(X_t,X_{t+1}) \in E_i \times H$, while positive for other possible $X_{t+1}$'s. Namely, for the $H$-steps (i.e., for $X_{t+1} \in \gamma_i \times H$), the Bayesian factors are all $1$, since $\cB_z$ depends only on the $\T^k$-coordinate, so its probability is the same starting from any vertex of $\gamma_i \times H$. For those $\T^k$-steps that are away from $\gamma$ and not excluded by $E_i$, the Bayesian factors are again $1$, because we will return to $\gamma_i \times H$ before making a step to $\gamma_{i\pm 1}\times H$. Finally, for the $\T^k$-step to $\gamma_{i-1}$ the Bayesian factor is at most $k$, by rearranging  
\begin{align*}
\Pb{\cB_z \md X_0,\dots,X_t=(\gamma_i,h) } 
&\ge \frac{1}{d+k}\Pb{\cB_z \md X_0,\dots,X_t=(\gamma_i,h),X_{t+1}=(\gamma_{i-1},h)}\\
&\qquad  + \frac{d}{d+k}\Pb{\cB_z \md X_0,\dots,X_t=(\gamma_i,h),X_{t+1} \in \gamma_{i} \times H)},
\end{align*}
which holds for any $t$ and any $h\in H$. Thus, the total weight of $H$-steps is $d$, while the total weight of $\T^k$-steps is at most $2k-1$, hence, before each  $\T^k$-step, the number of $H$-steps stochastically dominates a $\mathsf{Geom}\big((2k-1)/(2k-1+d)\big)-1$ variable. The stopping times $\alpha_i$ and $\beta_i$ are measurable with respect to the $\T^k$-coordinate of $(X_t)$, hence conditioned on all the events of~(\ref{e.backT}), i.e., on $\beta^\T_{i-1}-\beta^\T_i > 2 k^6$ and $\cB_z$ and on  $\cF(\beta_i)$, the variable $\beta^H_{i-1}-\beta^H_i $ stochastically dominates a sum of $2k^6$ iid variables with mean $d/(2k-1)$ and variance $d(2k-1+d)/(2k-1)^2$. Since $d\ge 2$, the expectation of the sum is larger than $2k^5$, and if $k$ is large enough, then the variance of the sum is less than $k^6$, hence the sum itself is larger than $k^5$ with a uniformly positive probability by Chebyshev's inequality.

For a proof of (\ref{e.backT}), first notice that, given $\cB_z$ and $\cF(\beta_i)$, the Bayesian factors calculated in the previous paragraph show that with a uniformly positive probability the step $(X_{\beta_i},X_{\beta_i+1})$ is in the $\T^k$-coordinate, away from $o$, into a branch different from $\gamma$ and $E_i$. (We are conditioning on the event $\cL_z$ of~(\ref{e.AL}) exactly in order for this uniformity to hold: the total Bayesian weight of these steps is at least $k/2-1$, while the total weight of all other steps is at most $k+d$.) From here, the distance of $(Y_t)$ from $\gamma_i$ is a biased random walk: whenever it changes (the step is in the $\T^k$-coordinate), it decreases with probability $1/k$ and increases otherwise. So, it will reach level $S_n$ with a uniformly positive probability. After this, whenever the walk is at $S_{n-1}$, it reaches level $S_i$ before $S_n$ only with probability $\asymp (k-1)^{i-n+1}$, by the usual exponential martingale argument \cite[Theorem 5.7.7]{Durrett}. For $i\leq n-7$, this is at most $O(k^{-6})$. That is, the number of steps in the $\T^k$-coordinate until returning to $S_i$ from $S_{n-1}$ stochastically dominates a geometric random variable with success probability $\asymp k^{-6}$, and this is at least $2k^6$ with a uniformly positive probability. This gives~(\ref{e.backT}).
\qed

The second ingredient will be that, both in $\pi_\mathrm{there}$ and $\pi_\mathrm{back}$, the amount of time spent in each bag $\gamma_i \times H$ is probably not very large. Namely, let $A_i$ be the set of times until time $\alpha_i$ when $X_t$ is in  $\gamma_i \times H$, and let $B_i$ be the analogous set of times from time $\beta_i$ until $\tau^+_o$. We let $H(A_i)$ and $H(B_i)$ be the set of vertices in $\gamma_i \times H$ visited at these times. Conditioned on $\cB_z$, the time spent in $\gamma_i \times H$, for $1 \leq i \leq n-1$, is stochastically dominated by a $\mathsf{Geom}\big(1/(k+d-1)\big)$ variable on the way to $z$ and by an independent copy on the way back to $o$, since the forward move along $\gamma$ is always available, the backward move is never,  and the forward move always has the largest Bayesian factor from~(\ref{e.Bayes}). The time spent in $\gamma_n \times H$ is $\mathsf{Geom}\big(1/(d+1)\big)$. So, letting $\cG_i$ denote the sigma-algebra generated by all the trajectory pieces {\it outside} the subgraph $G_i$ spanned by $\gamma_i \times H$ and the subgraphs of $G\setminus (\gamma\times H)$ hanging from there, up to time-translations for each piece (so, without the information how many steps within $G_i$ are taken), we have that, for any small $\delta>0$, if $D$ is a large enough absolute constant, then, for $i=1,\dots,n$,

\be\label{e.guests}
\begin{aligned}
\PB{|A_i|, |B_{i} | <  Dk \md  \cB_z, \cG_i } 
&\ge  \PB{\mathsf{Geom}\big(1/(k+d-1)\big) < Dk}^2\\
&= \left(1-\left(1-\frac{1}{k+d-1}\right)^{Dk}\right)^2\\
& > 1-\delta\,.
\end{aligned}
\ee

Now, if we have $|A_{i-1}|, |B_{i-1} | <  Dk$ and also the event $\{\beta^H_{i-1}-\beta^H_i > k^5 \}$ of~Lemma~\ref{l.back}, then $t-s > k^5$ for all $s\in A_{i-1}$ and $t\in B_{i-1}$, and hence the following lemma will be relevant to achieving $H(A_{i-1})\cap H(B_{i-1})=\emptyset$.

\begin{lemma}\label{l.heatkernel}
In any $d$-regular finite graph $H$ on more than $k^{5/2}$ vertices, if $t>k^5$, and $x,y\in V(H)$ are arbitrary, then the simple random walk heat kernel satisfies $\Ps{X_t=y \md X_0=x} < C_d k^{-5/2}$, with a constant $C_d<\infty$ that depends only on $d$.
\end{lemma}

\proof This is basically a special case of \cite[Lemma 3.6 in the arXiv version]{Russ:enum} or \cite{evolving}, with a few minor additional remarks.

In both references, the Markov chain is supposed to have a uniform laziness. So, we apply these results to the chain given by two consecutive steps on $H$. Since $H$ is $d$-regular, the probability of staying put in this chain is $1/d$. The stationary distribution is uniform. So, the references imply the on-diagonal bound $\Ps{X_{2t}=x \md X_0=x} < C_d k^{-5/2}$ for all even times $2t > k^{5/2}$. We then get the same off-diagonal bound $\Ps{X_{2t}=y \md X_0=x} < C_d k^{-5/2}$ by a standard Cauchy-Schwarz argument and the uniformity of the stationary distribution. Finally, to get the same bound for $X_{2t+1}$ being at $y$, average the bound over the neighbors of $y$ at time $2t$, before making the last step.

We also remark that to apply \cite{evolving} one has to take $\eps=k^{-5/2}|H|$ there, which is not small (as suggested by the notation $\eps$), but that is actually not a requirement in that paper.
\qed

If the trajectory $(X_t)_{t=0}^{\tau^+_o}$ satisfies $\cB_z$ and the intersection $\bigcap_{i=2}^{n-7} \big\{ H(A_{i-1}) \cap H(B_{i-1}) = \emptyset \big\}$, then its loop-erasure will intersect $\gamma_{n-7} \times H$, implying the event that we are interested in. We will give an exponentially small lower bound on the probability of this event, with a base that does not depend on $k$.

First of all, let $\Good_{n-7}:=\bigcap_{i=1}^{n-8} \big\{ |A_i | < Dk \big\} \in \cF(\alpha_{n-8})$, and then iteratively, for $i\leq n-7$,
\begin{align*}
\Prep_{i-1} &:= \Good_i \cap \big\{ \beta^H_{i-1}-\beta^H_i > k^5 \big\} \in \cF(\beta_{i-1})\,,\\ 
\Good_{i-1} &:= \Prep_{i-1}  \cap \big\{ H(A_{i-1}) \cap H(B_{i-1}) = \emptyset \big\} \in \cF(\beta_{i-2})\,.
\end{align*}
By~(\ref{e.guests}), we have
\begin{equation}\label{e.n-7}
\Pb{\Good_{n-7} \md \cB_z} > (1-\delta)^{n-8}. 
\end{equation}
We will now give a lower bound on $\Pb{\Good_{i-1} \md \cB_z, \, \Good_i}$, for each $2\leq i \le n-7$.
\medskip

First, consider $i=n-7$. By Lemma~\ref{l.back}, we have $\Pb{ \Prep_{n-8} \md \cB_z,\, \cF(\beta_{n-7})} > b$. Conditioned on $\cB_z$ and $\cF(\beta_{n-8}) \cap \Prep_{n-8}$, the bound~(\ref{e.guests}) gives that $|B_{n-8}| < Dk$ also holds with probability at least $1-\delta$. Finally, since the actual $H$-steps taken in the walk $(X_t)_{t\ge 0}$ are independent of the $\T^k$-steps and of the number of $H$-steps, Lemma~\ref{l.heatkernel} gives 
$$
\PB{ H(A_{n-8}) \cap H(B_{n-8}) \not= \emptyset \md  \cB_z,\, \cF(\beta_{n-8}) ,\, \Prep_{n-8},\, \big\{|B_{n-1} |  <  Dk\big\} } \\
< (Dk )^2 C_d k^{-5/2} < \delta,
$$
for our earlier small $\delta>0$, provided that $k$ is large enough. Altogether, we have
\begin{equation}\label{e.n-8}
\PB{ \Good_{n-8} \md  \cB_z,\, \cF(\beta_{n-8}) ,\,  \Good_{n-7} } >  (1-\delta)^2 b,
\end{equation}
finishing the first step of the induction. Also, note for future reference that we have proved, for every $i \leq n-8$, that
\begin{equation}\label{e.prepgood}
\PB{ \Good_{i} \md  \cB_z,\, \cF(\beta_{i}),\, \Prep_{i} } >  (1-\delta)^2 > 1-2\delta.
\end{equation}
%

The complication for the general $i \le n-8$ step to get from $\Good_i$ to $\Good_{i-1}$ will be that the event $\big\{ H(A_{i-1}) \cap H(B_{i-1}) = \emptyset \big\} \supset \Good_{i-1}$ concerns the trajectory from $\beta_{i-1}$ till $\beta_{i-2}$, hence interferes with the ``future'' event $\big\{ \beta^H_{i-2}-\beta^H_{i-1} > k^5 \big\} \supset \Prep_{i-2}$.
\medskip

Now fix any $i \le n-8$. We use Lemma~\ref{l.back} to get 
\begin{align}
\Pb{\Prep_{i-1} \md \cB_z ,\, \cF(\beta_i) ,\, \Good_{i}} 
& > \Pb{ \Prep_{i-1} \md \cB_z ,\, \cF(\beta_i) ,\, \Prep_i } -  \Pb{ \Good_{i}^c \md \cB_z ,\, \cF(\beta_i) ,\, \Prep_i } \nonumber \\
& > b - 2\delta, \label{e.goodbeta}
\end{align}
where the upper bound $2\delta$ is from~(\ref{e.prepgood}). If we set $\delta:=b/4 \leq 1/4$, then this lower bound becomes $b/2 > 0$, and it makes sense to continue as follows:
\begin{align}
\Pb{\Good_{i-1} \md \cB_z ,\, \cF(\beta_i) ,\, \Good_{i},\, \Prep_{i-1} } \nonumber \\
& \hskip -1 in > \Pb{ \Good_{i-1} \md \cB_z ,\, \cF(\beta_i) ,\, \Prep_{i-1} } -  \Pb{ \Good_{i}^c \md \cB_z ,\, \cF(\beta_i) ,\, \Prep_{i-1} } \nonumber \\
& \hskip -1 in > (1-\delta)^2 - \frac{  \Pb{ \Good_{i}^c \md \cB_z ,\, \cF(\beta_i) ,\, \Prep_{i} } }{  \Pb{ \Prep_{i-1} \md \cB_z ,\, \cF(\beta_i) ,\, \Prep_{i} } } \nonumber \\
& \hskip -1 in \ge (1-\delta)^2 - \frac{2\delta}{b} \ge \frac{9}{16}-\frac{1}{2} =\frac{1}{16}.
\label{e.goodgood}
\end{align}
The estimates~(\ref{e.goodbeta}) and~(\ref{e.goodgood}) together give
\be\label{e.bg}
\Pb{ \Good_{i-1} \md  \cB_z,\, \cF(\beta_i),\,\Good_i } > \frac{b}{32},
\ee
for $2 \leq i \leq n-8$. Telescoping this with~(\ref{e.n-7}) and~(\ref{e.n-8}), we get
\be\label{e.bgtele}
\PB{\bigcap_{i=2}^{n-7} \big\{ H(A_{i-1}) \cap H(B_{i-1}) = \emptyset \big\} \md \cB_z } \ge \PB{ \bigcap_{i=1}^{n-7} \Good_i \md \cB_z} \ge  \left(\frac{3}{4}\right)^{n-8} \frac{9b}{16} \left(\frac{b}{32}\right)^{n-9}.
\ee

From this exponentially small lower bound, to find a good $z \in S_n$ and hence a good $\tilde z \in S_{n-7}$ with a uniformly positive probability,  we will again use the second moment method, for which we need a little bit of preparation. Define the events 
\begin{align*}
\cC_z &:= \cB_z \cap \big\{ H(A_i) \cap H(B_i) = \emptyset \textrm{ for all } i=1,\dots, n-8 \big\},\\
\cC_z(h) &:= \cC_z \cap \{ X_{\tau^+_o}= (o,h) \},\quad h\in H.
\end{align*}
Then, as~(\ref{e.bgtele}) says,
\be\label{e.bgn}
\Pb{\cC_z \md \cB_z} > \frac{9}{16} \left( \frac{3b}{128}\right)^{n-8}.
\ee
Furthermore, we claim that
\be\label{e.minmax}
\max_{h\in H} \Pb{\cC_z(h) \md \cB_z} \leq C \min_{h\in H} \Pb{\cC_z(h) \md \cB_z},
\ee
with a constant $C<\infty$ that may depend on $H$ and $k$, but not on $n$. 

In the proof of this claim, we will use a small technical lemma:

\begin{lemma}\label{l.2conn}
Every finite transitive graph $H$ is 2-vertex-connected: for any vertex $g \in V(H)$, the graph we get from $H$ by deleting $g$ is still connected.
\end{lemma}

\proof Assume that there is a cut-vertex $g$, whose removal cuts $H$ into at least two components; denote the largest of these  by $H_g$ (or one of the largest ones in case of a draw). Take some vertex $h$ not in $\{g\} \cup H_g$. By transitivity, $h$ is also a cut-vertex, whose removal results in at least two components, one containing both $g$ and $H_g$. But this component will have a size strictly larger than $|H_g|$, contradicting transitivity.\qed

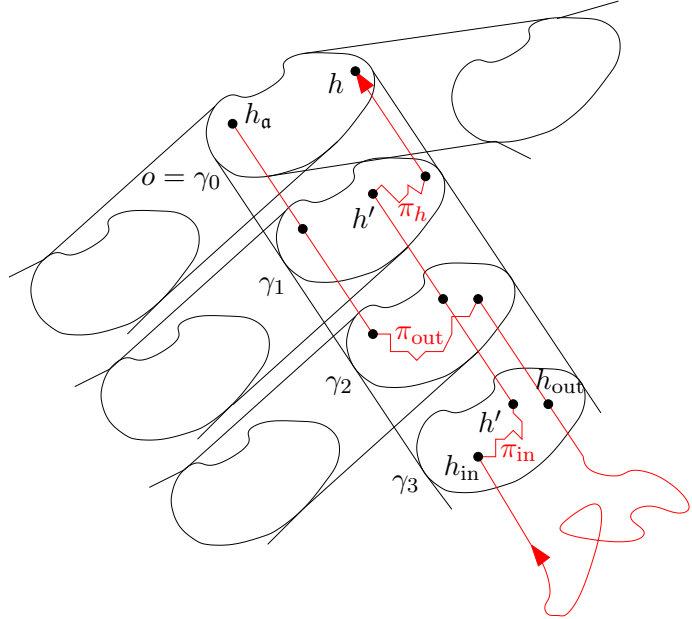
\begin{figure}[htbp]
\SetLabels
(0.26*0.7) $o=\gamma_0$\\
(0.39*0.53) $\gamma_1$\\
(0.485*0.37) $\gamma_2$\\
(0.58*0.21) $\gamma_3$\\
(0.37*0.8) $h_{\fa}$\\
(0.48*0.85) $h$\\
(0.52*0.63) $h'$\\
(0.7*0.3) $h'$\\
(0.66*0.23) $h_{\textrm{in}}$\\
(0.8*0.37) $h_{\textrm{out}}$\\
(0.74*0.26) $\textcolor{red}{\pi_{\textrm{in}}}$\\
(0.60*0.45) $\textcolor{red}{\pi_{\textrm{out}}}$\\
(0.59*0.65) $\textcolor{red}{\pi_h}$\\
\endSetLabels
\centerline{
\AffixLabels{
\includegraphics[width=0.55\textwidth]{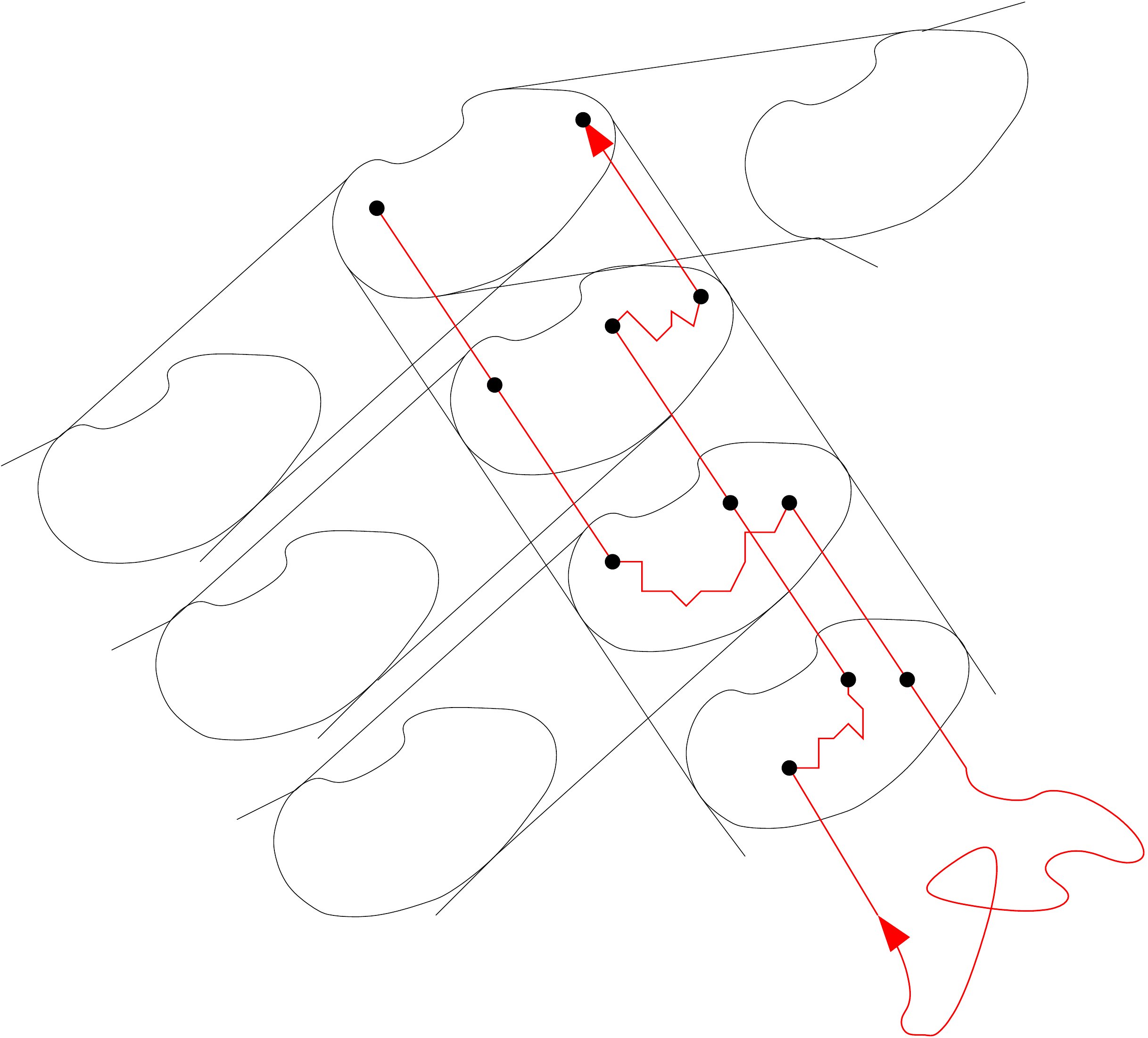}
}}
\caption{Producing a good random walk trajectory in $\T^k \times H$.}
\label{f.rewiring}
\end{figure}

To prove~(\ref{e.minmax}), first observe that, if we condition the random walk trajectory $(X_t)$ to satisfy $\cC_z$, and let $X_{\alpha_3}=(\gamma_3,h_{\mathrm{out}})$ be the last vertex in $\gamma_3\times H$ on the trajectory before $\tau_z$, and let $X_{\beta_3}=(\gamma_3,h_{\mathrm{in}})$ be the first one after $\tau_z$, then, conditionally on $h_{\mathrm{out}}$ and $h_{\mathrm{in}}$, the part of the trajectory between $h_{\mathrm{out}}$ and $h_{\mathrm{in}}$ is independent of the rest. Therefore, if we prove that there exists some $p>0$, depending only on $H$ and $k$, but not on $n$, such that, for any two vertices $h_{\mathrm{out}} \not= h_{\mathrm{in}}$, and any $h\in H$, the probability that $(X_t)_{t=0}^{\alpha_3}$ and $(X_t)_{t=\beta_3}^{\tau^+_o}$ satisfy the conditions of $\cC_z(h)$ relating to $\gamma_i\times H$ for $i=0,1,2,3$ is at least $p$, then $C=1/p$ will clearly work in~(\ref{e.minmax}). For the argument that follows, see Figure~\ref{f.rewiring}. 

Pick any vertex $h' \in H\setminus\{h_{\mathrm{out}},h_{\fa} \}$. By the 2-connectedness of $H$, we can pick a path $\pi_{\mathrm{in}}$ in $\gamma_3\times H$ between $(\gamma_3,h_{\mathrm{in}})$ and $(\gamma_3,h')$ that avoids $(\gamma_3,h_{\mathrm{out}})$, a path $\pi_{\mathrm{out}}$ in $\gamma_2\times H$ between $(\gamma_2,h_{\mathrm{out}})$ and $(\gamma_2,h_\fa)$ that avoids $(\gamma_2,h')$, and a path $\pi_h$ in $\gamma_1\times H$ between $(\gamma_1,h)$ and $(\gamma_1,h')$ that avoids $(\gamma_1,h_\fa)$. Then $(X_t)_{t=0}^{\alpha_3}$ can go from $(\gamma_0,h_\fa)$ straight to $(\gamma_2,h_\fa)$, then to $(\gamma_2,h_{\mathrm{out}})$ via  $\pi_{\mathrm{out}}$, then straight to $(\gamma_3,h_{\mathrm{out}})$, and $(X_t)_{t=\beta_3}^{\tau^+_o}$ can go from $(\gamma_3,h_{\mathrm{in}})$ via $\pi_{\mathrm{in}}$ and $\pi_h$ to $(\gamma_0,h)$. All of this happens with probability at least $(d+k)^{-3|H|-6}$, which proves~(\ref{e.minmax}). (Note that we needed the extra vertex $h'$ and the four layers $\gamma_0,\dots,\gamma_3$ for this construction because it might happen that
 $h=h_{\mathrm{out}}$; otherwise, taking $h':=h$ and removing the $\gamma_1$ layer could have worked.)
\medskip

We are now ready for the second moment method. Let $W_n$ be the set of leaves $z\in S_n$ that satisfy $\cC_z(h_\fb)$, with the desired endpoint $\fb=(o,h_\fb)$. We will run a second moment argument, as in Section~\ref{s.tree}, to show that $W_n$ is non-empty with a positive probability, uniformly in $n$. 

First note that~(\ref{e.bgn}) and~(\ref{e.minmax}) imply
\be\label{e.qn}
q(n):=\Pb{ \cC_z(h_\fb) \md \cB_z} > c \, \left(\frac{3b}{128}\right)^n, 
\ee
where $c$ depends on $H$ and $k$, but not on $n$. Together with~(\ref{e.Bz}), we have
\be\label{e.W1}
\E |W_n| \, \asymp_{k,H} \, (k-1)^n\, p_k^{n}\, q(n).
\ee
Using~(\ref{e.pk}) and~(\ref{e.qn}), this tends to infinity as $n\to\infty$ for $k$ large enough. 

To estimate the second moment $\Eb{|W_n|^2}$, let $z,v \in S_n$ be leafs such that their last common ancestor is $w\in S_m$, with $m\ge 1$. We claim that
\be\label{e.qq}
\Pb{\cC_z(h_\fb) \cap \cC_v(h_\fb) \md \cB_z \cap \cB_v} \leq Q \, q(n)\, q(n-m),
\ee
with some $Q<\infty$ that depends only on $H$ and $k$, but not on $n$. 

First note that we may assume that $m\leq n-10$, since otherwise the factor $q(n-m)$ on the right hand side of~(\ref{e.qq}) is obviously at least a positive constant that does not depend on $n$, hence a suitable $Q$ does exist.

By symmetry, we may assume $\tau_z < \tau_v$. We first show that 
\be\label{e.coupling}
\Pb{\cC_z(h_\fb) \md \cB_z \cap \cB_v \cap \{\tau_z < \tau_v\} } \leq C_{k,H} \, \Pb{\cC_z(h_\fb) \md \cB_z }.
\ee
We do this by coupling (with a uniformly positive probability) the trajectory $(X_t)$ conditioned on $\cB_z$ to be identical to the trajectory conditioned on $\cB_z \cap \cB_v \cap \{\tau_z < \tau_v\}$, denoted by $(\tilde X_t)$, within the ray $\gamma \times H$ (which leads from $o$ to $z$), except for a bounded neighborhood of $w=\gamma_m$. Given $\cB_z$, we know from~(\ref{e.bg}) that $H(A_m)\cap H(B_m)=\emptyset$ occurs with a uniformly positive probability, say $c_1>0$. Conditioning on  $\cB_z \cap \{H(A_m)\cap H(B_m)=\emptyset\}$ gives a certain distribution to the pairs of vertices $\big(X_{\alpha_{m-1}},X_{\beta_{m-1}}\big)$ and $\big(X_{\alpha_{m}},X_{\beta_{m}}\big)$, which are basically the places where the trajectory leaves $w\times H$.  
%
On the other hand, conditioning on $\cB_z \cap \cB_v \cap \{\tau_z < \tau_v\}$, we get some distribution on 
$(\tilde X_t)_{t=\alpha_{m-4}}^{\alpha_{m+3}}$ and  $(\tilde X_t)_{t=\beta_{m+3}}^{\beta_{m-4}}$.  
Whenever these pieces of $(\tilde X_t)$-trajectories satisfy $H(A_i)\cap H(B_i)=\emptyset$ for $i=m-3,\dots,m+3$ (so that $(\tilde X_t)$ still has a chance to satisfy $\cC_z(h_\fb)$), the argument of Figure~\ref{f.rewiring} gives that, conditioned on these trajectory pieces, with a probability at least $c_2>0$ that depends only on $k$ and $H$, we have that $(X_t)_{t=\alpha_{m-4}}^{\alpha_{m+3}}$ and  $(X_t)_{t=\beta_{m+3}}^{\beta_{m-4}}$ satisfy $\big(X_{\alpha_{m + 3}}, X_{\beta_{m+ 3 }}\big)=\big( \tilde X_{\alpha_{m + 3}}, \tilde X_{\beta_{m+ 3 }}\big)$ and $\big(  X_{\alpha_{m - 4}},  X_{\beta_{m- 4}}\big)= \big( \tilde X_{\alpha_{m - 4}}, \tilde X_{\beta_{m- 4}}\big)$. Conditioned on these equalities, we can couple  the trajectories $(X_t)_{t=0}^{\alpha_{m-4}}$, $(X_t)_{t=\alpha_{m+3}}^{\beta_{m+3}}$, and $(X_t)_{t=\beta_{m-4}}^{\tau_o^+}$ to be equal to the tilde versions, hence if $(\tilde X_t)$ satisfies $\cC_z(h_\fb)$, so does $(X_t)$. Altogether,~(\ref{e.coupling}) follows with $C_{k,H}=1/(c_1c_2)$.
 %

Now let $H(A_m), H(B'_m), H(B''_m)$ be the set of vertices in $w\times H$ visited before $\tau_z$, between $\tau_z$ and $\tau_v$, and after $\tau_v$, respectively; thus $H(B'_m) \cup H(B''_m) = H(B_m)$. Notice that $\cC_z(h_\fb) \cap \cC_v(h_\fb)$ implies that $H(A_m), H(B'_m), H(B''_m)$ are mutually disjoint, an event we will denote by $\cM_w$. Condition now, beyond $\cB_z \cap \cB_v \cap \{\tau_z < \tau_v\}$, also on the event $\cC_z(h_\fb) \cap \cM_w$. Let $h'$ be the vertex in $H(B'_m)$ last visited before $\tau_v$, and let $h''$ be the first vertex in $H(B''_m)$ after $\tau_v$. Since $H$ is transitive, there is an automorphism taking $h'$ to $h_\fa$, and $h''$ to some $h^*$. Now, the events along the ray from $w$ to $v$ that are needed for $\cC_v(h_\fb)$ are just the events for some length $n-m$ ray, with the extra condition that the first step from $(w,h_\fa)$ and the last step to $(w,h^*)$ are both in the $\T^k$-coordinate. Thus, using~(\ref{e.minmax}), we have 
\be\label{e.Mw}
\PB{\cC_v(h_\fb) \md \cC_z(h_\fb) \cap \cM_w \cap \cB_z \cap \cB_v \cap \{\tau_z < \tau_v\} } < C'_{k,H}\, q(n-m).
\ee
Since $\cM_w \supset \cC_z(h_\fb) \cap \cC_v(h_\fb)$, we can combine this with
$$
\Pb{\cC_z(h_\fb) \cap \cM_w \md \cB_z \cap \cB_v \cap \{\tau_z < \tau_v\} } \leq C_{k,H} \, q(n),
$$
which we get from~(\ref{e.coupling}), and we arrive at~(\ref{e.qq}).

From~(\ref{e.qq}) and~(\ref{e.BzBv}), similarly to~(\ref{e.Z2}), we have
\be\label{e.W2}
\E |W_n|^2 \leq Q' \sum_{m=1}^n k(k-1)^{m-1}\, (k-1)^{2(n-m)}\, p_k^{2n-m} \, q(n)\,q(n-m).
\ee
For the Cauchy-Schwarz second moment method, we want that $\E |W_n|^2 < Q'' (\E |W_n|)^2$, for some $Q''<\infty$ that does not depend on $n$. Substituting~(\ref{e.W1}) and~(\ref{e.W2}) into this inequality, then rearranging, we arrive at the following inequality to prove: 
\be\label{e.Q3}
\sum_{m=1}^n \big((k-1)p_k\big)^{-m} \,q(n-m) \stackrel{?}{<} Q''' \, q(n).  
\ee
The final ingredient is that, writing $y$ for vertex $\gamma_{m}$ on the ray from $o$ to $z$, again only for $m\leq n-10$, and writing $\cB_z^y$ and  $\cC_z^y$ for the analogs of the events $\cB_z=\cB_z^o$ and  $\cC_z=\cC_z^o$ when the root is $y$ instead of $o$, 
\begin{align*}
\frac{q(n)}{q(n-m)} = \frac{ \Ps{\cC_z^o(h_\fb) \md \cB_z^o} }{ \Ps{\cC_z^y(h_\fb) \md \cB_z^y} } 
&\asymp_{k,H} \frac{ \Pb{ H(A_i) \cap H(B_i) = \emptyset \textrm{ for } i=1,2,\dots, n-8 \md \cB_z^o} }
{ \Pb{ H(A_i) \cap H(B_i) = \emptyset \textrm{ for } i=m+1,m+2,\dots, n-8 \md \cB_z^y}  }\\
&\asymp_{k,H} \frac{ \Pb{ H(A_i) \cap H(B_i) = \emptyset \textrm{ for } i=1,2,\dots, n-8 \md \cB_z^o} }
{ \Pb{ H(A_i) \cap H(B_i) = \emptyset \textrm{ for }  i=m+1,m+2,\dots, n-8 \md \cB_z^o}  }\\
&= \Pb{ \cC_z^o \md \cB_z^o,\ H(A_i) \cap H(B_i) = \emptyset \textrm{ for } i=m+1,m+2,\dots, n-8} \\
& > (3b/128)^m, 
\end{align*}
where the first $\asymp$ is by~(\ref{e.minmax}); the second $\asymp$ is by a coupling argument similar to the one that gave~(\ref{e.coupling}), now doing the coupling in $\{\gamma_m,\dots,\gamma_{m+3}\}\times H$; and the inequality in the last line follows from~(\ref{e.n-7}) and~(\ref{e.bg}), just like in~(\ref{e.bgtele}). Plugging this into~(\ref{e.Q3}), we arrive at
$$
\sum_{m=1}^n \big((k-1)p_k\big)^{-m} (2/bg)^m \stackrel{?}{<} Q'''',
$$
which is true if $k$ is large enough, since~(\ref{e.pk}) tells us that $(k-1)p_k\to\infty$ as $k\to\infty$. This finishes the proof of the disconnectedness Theorem~\ref{t.disco}. \qed

For the first direct proof of having infinitely many trees almost surely, pick an infinite ray $o_1,o_2,\dots$ in $\T^k$, pick any $h\in H$, and let $\fa_i:=(o_i,h)$. Our exhaustion $G_n=T_n\times H$ contains $\fa_1,\dots,\fa_n$. Perform Wilson's algorithm in $G_n$ as follows. 

First run a $\LERW$ from $\fa_2$ to $\fa_1$, denoted by $\ell_1$. By a small modification of our previous proof, with a positive probability that depends only on $H$ and $k$, this $\ell_1$ will first enter the subtree (times $H$) that starts at $o_2$ and does not contain $o_1$ or $o_3$, then will hit the boundary $S_n \times H$, then hits $o_2\times H$ at a vertex $\fb_2=(o_2,h_2)$ different from $\fa_2$, then goes straight to $(o_1,h_2)$, then hits $\fa_1=(o_1,h)$ without leaving $o_1\times H$. Without conditioning on this good event, denoted by $\cG_1$ hereafter, the $\T^k$-coordinate of the random walk that gives $\ell_1$,  
viewed only at the times when it moves on the ray $o_1,\dots,o_n$, performs a simple random walk on this segment until $\tau^+_o$. 
The maximum $j$ for which $o_j$ is touched by the projection is stochastically dominated by the maximum of a one-dimensional random walk excursion, which is almost surely finite, since the walk is recurrent.
The maximum $j$ for which $o_j\times H$ is touched by $\ell_1$, denoted by $j_2$, is even smaller. Let $\fb_2=(o_{j_2},h_2)$ be the last vertex in $j_2\times H$ touched by $\ell_1$. Note that this definition of $\fb_2$ extends our previous one that we made under $\cG_1$.

Next, run a $\LERW$ from $\fa_{j_2+1}$ to $\fb_2$, denoted by $\ell_2$, which, with a  positive probability that depends only on $H$ and $k$, will enter the subtree (times $H$) that starts at $o_{j_2+1}$ and does not contain $o_{j_2}$ or $o_{j_2+2}$, then will hit the boundary of $T_n \times H$, then hits $o_{j_2+1}\times H$ at a vertex $\fb_3=(o_{j_2+1},h_3)$ different from $\fa_{j_2+1}$, then goes straight to $(o_{j_2},h_3)$, then hits $\fb_2$ without leaving $o_{j_2}\times H$. Without conditioning on this good event, denoted by $\cG_2$, the maximum $j$ for which $o_j\times H$ is touched by $\ell_2$, denoted by $j_3$, has the property that $j_3-j_2$ is stochastically dominated by the maximum of a one-dimensional simple random walk excursion. Let $\fb_3=(o_{j_3},h_3)$ be the last vertex in $j_3\times H$ touched by $\ell_2$, extending the definition that we made under $\cG_2$.

Iterate this procedure until we have reached $o_n\times H$, producing the $\LERW$ paths $\ell_1,\dots,\ell_{I_n}$. Since the distribution of $j_{i+1}-j_i$ is always stochastically dominated by the maximum of a one-dimensional simple random walk excursion, the variable $I_n$ tends to infinity in probability, as $n\to\infty$. Each $\ell_i$, independently of the previous ones, satisfies $\cG_i$ with a  positive probability that depends only on $H$ and $k$. Thus, the number of events $\cG_i$ satisfied also tends to infinity in probability. This shows that the number of trees in the weak limit is almost surely infinite.

\section{You can't hide from yourself}\label{s.change}

\proofof{Theorem~\ref{t.genset}} 
The natural free generating set in each coordinate of the product, together with their inverses, gives a tree $\T^{2k}$ in the $\F_k$ coordinate and a cycle in the $H=\Z_{k^9}$ coordinate (every edge that appears does so in both orientations, so, as usual, we consider them to be unoriented single edges). If $k$ is large enough, then Theorem~\ref{t.disco} tells us that the $\FUSF$ has infinitely many trees almost surely.

The second Cayley graph will also be a direct product graph: we again take free generators for $\F_k$ with their inverses, while all the elements in $H=\Z_{k^9}$, except for the identity. This gives the Cayley graph $\T^{2k} \times K_{k^9}$, where $K_n$ is the complete graph on $n$ vertices with a single unoriented edge between any pair of vertices.

We will show that, for the $\LERW$ from ${\fa}=(o,h_{\fa})$ to ${\fb}=(o,h_{\fb})$, with $h_{\fa}\not=h_{\fb}\in H$, the probability that the $\LERW$ is not contained in $T_r \times H$ is exponentially small in $r$, if $k$ is large enough. (As before, $T_r$ is the ball of radius $r$ in $\T^{2k}$.) This of course implies the theorem.

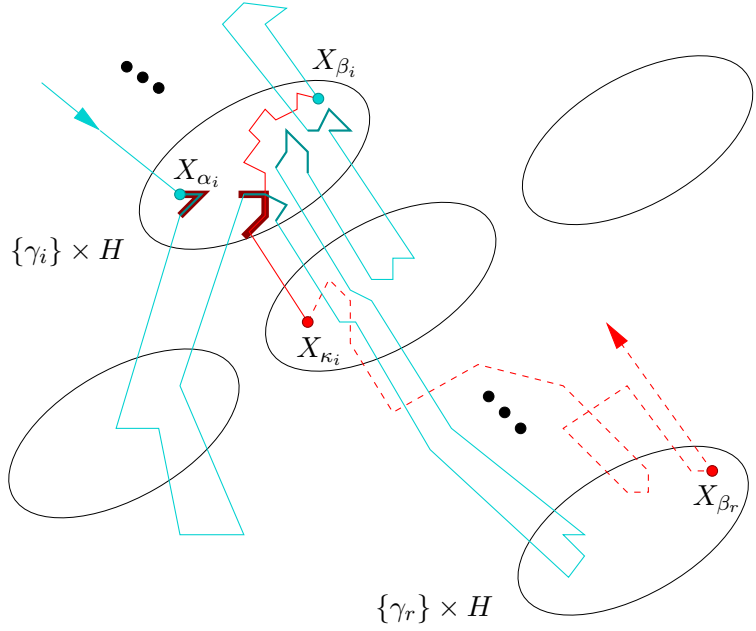
\begin{figure}[htbp]
\SetLabels
(0.09*0.58) $\{\gamma_i\} \times H$\\
(0.57*0) $\{\gamma_r\} \times H$\\
(0.26*0.71) $X_{\alpha_i}$\\
(0.44*0.89) $X_{\beta_i}$\\
(0.42*0.42) $X_{\kappa_i}$\\
(0.94*0.18) $X_{\beta_r}$\\
\endSetLabels
\centerline{
\AffixLabels{
\includegraphics[width=0.6\textwidth]{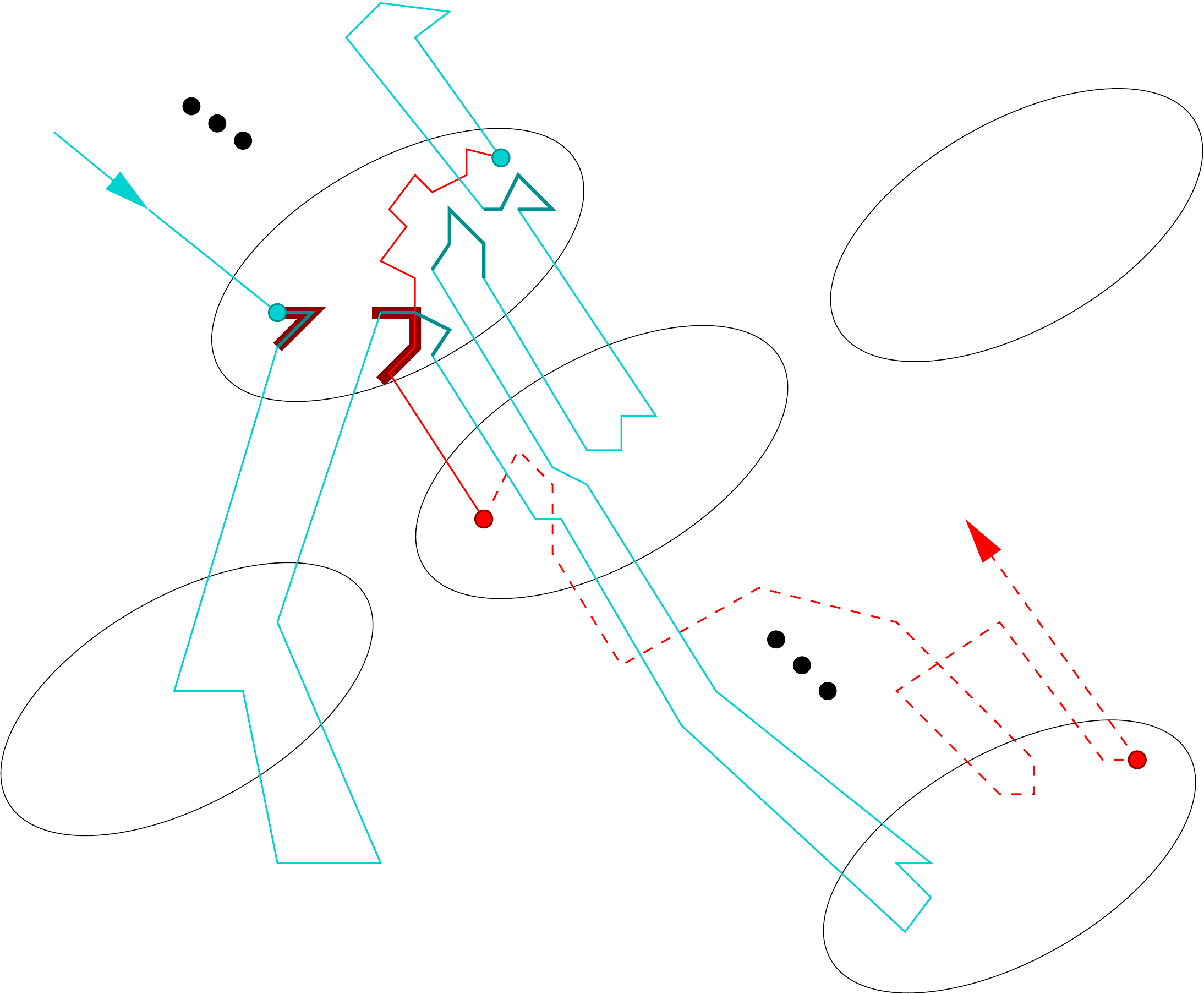}
}}
\caption{The cyan-colored path is $\LERW_{\beta_i}$. Its darker cyan parts are $\LERW_{\beta_i} \cap (\{\gamma_i\} \times H)$. The solid red path is $(X_t)_{t=\beta_i}^{\kappa_i}$, the dashed red path is $(X_t)_{t\ge \kappa_i}$  and the thick dark red pieces form $\LERW_{\kappa_i} \cap (\{\gamma_i\} \times H)$.}
\label{f.abkappa}
\end{figure}

Fix any ray $o=\gamma_0, \gamma_1,\dots,\gamma_r$ in $T_r$, and let $\beta_r$ be the last time that the simple random walk $(X_t)_{t=0}^{\tau_\fb}$  from ${\fa}$ to ${\fb}$ enters the bag $\{\gamma_r\} \times H$ (i.e., the last $\beta$ such that $X_{\beta}\in \{\gamma_r\} \times H$ but $X_{\beta-1}\not\in \{\gamma_r\} \times H$). If the walk never enters $\{\gamma_r\} \times H$, we set $\beta_r=\infty$. If $\beta_r<\infty$, then we also let $\beta_i$, for $i=0,1,\dots,r-1$, be the last time before $\beta_r$ that $(X_t)_{t\ge 0}$ enters the bag $\{\gamma_i\} \times H$, and let $\kappa_i$ be the first time after $\beta_i$ that $(X_t)_{t\ge 0}$ is not in $\{\gamma_i\} \times H$ (in fact, because of $\beta_r<\infty$, we have $X_{\kappa_i} \in \{\gamma_{i+1}\} \times H$). 
Furthermore, we let $\LERW_t$ denote the loop-erasure of $(X_s)_{s=0}^t$, and, still assuming $\beta_r<\infty$, we let $\alpha_i$ be the first time $\alpha$ with the property that $X_\alpha \in \{\gamma_i\} \times H$ and $\LERW_{\alpha} \cap \{X_s\}_{s=\alpha}^{\beta_i} = \{ X_\alpha \}$. In other words, $X_{\alpha_i}$ is the first vertex along $\LERW_{\beta_i}$ that is in $\{\gamma_i\} \times H$. See Figure~\ref{f.abkappa}.

We will show that, with very high probability conditionally on $\beta_r<\infty$, for most $i=0,1,\dots,r$ the intersection $\LERW_{\kappa_i} \cap (\{\gamma_i\} \times H)$ is quite large. This will imply, introducing the notation $\varphi_i$ for the first time when the simple random walk from $\beta_r$ to $\tau_{\fb}$ enters $\{\gamma_i\} \times H$ again, and $\psi_i$ for the first time after $\varphi_i$ when the simple random walk is in $\{\gamma_{i\pm 1}\} \times H$, that the event 
\be\label{e.avoidone}
\Avoid_i := \{\beta_r < \infty\} \cap \left\{  \LERW_{\kappa_i} \cap (X_t)_{t=\varphi_i}^{\psi_i} = \emptyset\right\}
\ee
is very unlikely to happen; in fact, we will show that
\be\label{e.avoidall}
\PB{\bigcap_{i=1}^{r-1} \Avoid_i} \leq \left(\frac{1}{2k}\right)^{r-1}.
\ee
This is relevant because any $\Avoid_i^c$ implies that the $\LERW$ from ${\fa}$ to ${\fb}$ does not intersect $\{\gamma_r\} \times H$. Thus, in order for the $\LERW$ from $\fa$ to $\fb$ not to be contained in $T_r\times H$, there must exist a ray $\gamma_0, \gamma_1,\dots,\gamma_r$ so that the event of~(\ref{e.avoidall}) occurs. Since the number of possible such rays is $(2k)(2k-1)^{r-1}$, a union bound using~(\ref{e.avoidall}) gives an exponentially small upper bound $2k\big(1-1/(2k)\big)^{r-1}$, as desired.

%

In the proof of~(\ref{e.avoidall}),  we will be conditioning on $\beta_r<\infty$ from now on (which only raises the probability). As hinted above, the key statement will be that $\LERW_{\kappa_i} \cap (\{\gamma_i\} \times H)$ is large with high probability, independently of other bags. More precisely, denoting by $\calf_i$ the sigma-algebra generated by $(X_t)_{t=0}^{\kappa_i}$ (including the value of $\kappa_i$), we will prove that 
\be\label{e.large}
\PB{ \big|\LERW_{\kappa_i} \cap (\{\gamma_i\} \times H)\big| < k^3 \md \calf_{i-1},\ \beta_r<\infty} < O(k^{-3}),
\ee
for any $i \in \{0,1,\dots,r-1\}$. For this, we will need a small Markov chain mixing time lemma that controls the size of the $\LERW$ in a complete graph. For basic definitions, such as the total variation distance $d_\mathrm{TV}$, see \cite{LPW}. 

\begin{lemma}\label{l.mixing}
Let $(X_t)_{t\ge 0}$ be simple random walk on the complete graph $K^{\circ}_n$ with loops; that is, each step of the walk is just a new independent vertex distributed as $\mathsf{Unif}\{1,\dots,n\}$. Now let $(L_t)_{t\ge 0}$ be the Markov chain on $\{1,\dots,n\}$ where $L_t$ is the size of the loop-erased version of the path $(X_s)_{s\ge 0}^t$. Then the following are true.
\begin{itemize}
\item[{\bf (1)}] The transition probabilities for $(L_t)_{t\ge 0}$ are $p(i,j)=1/n$ for all $j \in \{1,\dots,i\}$, and $p(i,i+1)=(n-i)/n$. The unique stationary distribution of the chain satisfies $\pi(i) \leq i/n$.
\item[{\bf (2)}] The total variation mixing time of the chain is $O(\sqrt{n})$; in fact, $d_\mathrm{TV}(\mu_t,\pi) < \exp\left(-\frac{t^2}{2n}\right)$ for every $t$, where $\mu_t$ is the distribution of $L_t$ started from any given state.
\item[{\bf (3)}] If $(X_t)_{t\ge 0}$ is simple random walk on the complete graph $K_n$ without loops, then the stationary distribution for the loop-erased version $L_t$ is the same as before, the mixing time is still $O(\sqrt{n})$, and hence $d_\mathrm{TV}(\mu_t,\pi) < \exp\left(-c\frac{t}{\sqrt{n}}\right)$ for some $c>0$, for every $t$ and every starting state.
\end{itemize}
\end{lemma}

\proof
{\bf (1)} At time $t$, if the next step $X_{t+1}$ is to the $j$th vertex on the current loop-erased path, then $L_{t+1}=j$; if $X_{t+1}$ is to a vertex not currently on the path, then $L_{t+1}=L_t+1$. The transition probabilities follow. This chain is clearly irreducible and aperiodic, hence it has a unique stationary distribution $\pi$, which satisfies the equation 
$$
\pi(i+1) = \frac{n-i}{n} \pi(i) + \frac{1}{n} \sum_{k=i+1}^n \pi(k) \leq \pi(i) + \frac{1}{n}.
$$
The inequality  $\pi(i) \leq i/n$ follows by induction on $i$.

{\bf (2)} We will bound the mixing time by a standard coupling argument: if $(L_t,\tilde L_t)_{t\ge 0}$ is any coupling of two copies of the Markov chain, one with $L_0=i$, the other with $\tilde L_0=j$, and $\tau_\mathrm{coupling}$ is the first time when $L_t=\tilde L_t$, then \cite[Corollary 5.5]{LPW} says that
\be\label{e.TVtau}
d_\mathrm{TV}(\mu_t,\pi) \leq \max_{i,j\in K^\circ_n} \Pb{\tau_\mathrm{coupling} > t}.
\ee
Our coupling will be a monotone one: we assume $i<j$, then will maintain $L_t \leq \tilde L_t$ for all $t\ge 0$. Take i.i.d.~random variables $U_t \sim \mathsf{Unif}\{1,\dots,n\}$ for $t > 0$. Given already $(L_s,\tilde L_s)_{s=0}^t$, we generate $(L_{t+1},\tilde L_{t+1})$ as follows. If $U_{t+1} \leq L_t$, then let $L_{t+1}:=U_{t+1}$; if $U_{t+1} > L_t$, then let $L_{t+1}:=L_t+1$. We make exactly the same definitions for $\tilde L_{t+1}$, using the same variable $U_{t+1}$ as for $L_{t+1}$. This is clearly a monotone coupling of two copies of the chain, and it has the property that $\tau_\mathrm{coupling}=\inf\{t+1: U_{t+1} \leq L_t+1\}$. (If $U_{t+1} \leq L_t < \tilde L_t$, then both chains are in the first case of the definition; if $L_t+1=U_{t+1} \leq \tilde L_t$, then $L_{t+1}$ is in the second case, $\tilde L_{t+1}$ is in the first case, but nevertheless they have become equal.) Therefore,
$$
 \Pb{\tau_\mathrm{coupling} > t} = \prod_{s=1}^{t}  \left(1-\frac{i+s}{n} \right) < \prod_{s=1}^{t} \left(1-\frac{s}{n} \right) < \exp \left(-\sum_{s=1}^{t} \frac{s}{n}\right) < \exp\left(-\frac{t^2}{2n}\right).
$$
This is true for any pair of starting states $1\leq i < j \leq n$, hence the result follows by~(\ref{e.TVtau}).

{\bf (3)} The non-laziness of the random walk $(X_t)$ causes only that the chain $(L_t)$ can never stay put; the new transition probabilities are the same as if we condition the old chain to actually move. Thus the stationary distribution remains the same. 

For the mixing time, first note that the TV-distance between $(L_{t+1} \md L_t=i)$ in the new chain and $(L^\textrm{old}_{t+1} \md L^\textrm{old}_t=i)$ in the old chain is $1/n$ for any $i$, hence we can couple the new chains $(L_t)$ and $(\tilde L_t)$ to the old chains $(L_t^\textrm{old})$ and $(\tilde L_t^\textrm{old})$ such that 
$$
\Pb{ L_{t+1} \not =L_{t+1}^\textrm{old} \md L_{t} =L_t^\textrm{old} } \leq \frac{1}{n},
$$
and similarly for the tilde versions. We now start the old and new chains at the same places $L_0=L_0^\textrm{old}=i$ and  $\tilde L_0=\tilde L_0^\textrm{old}=j$. Applying our previous coupling between $(L_t^\textrm{old})$ and $(\tilde L_t^\textrm{old})$ and the bound from part (2), we get that, for $t_0:=\lceil 2\sqrt{n} \rceil$,
\begin{align*}
\Pb{L_{t_0} \not= \tilde L_{t_0}} 
&\leq \Pb{L^\textrm{old}_{t_0} \not= \tilde L^\textrm{old}_{t_0}}  + 
\Pb{ L_i \not= L_i^\textrm{old} \textrm{ or }   \tilde L_i \not= \tilde L_i^\textrm{old} \textrm{ for some }i=1,\dots, t_0 }\\   
&\leq  \exp(-2) + \frac{2t_0}{n}.
\end{align*}
This is smaller than $1/4$ if $n$ is large enough, hence the mixing time is at most $t_0$, and the exponential decay of the TV-distance follows from a standard argument \cite[Section 4.5]{LPW}.
\qed

To prove~(\ref{e.large}, we do not only condition on $\calf_{i-1}$ and $\beta_r<\infty$, but also on $(X_t)_{t=0} ^{\beta_i}$, including the value of $\beta_i$. Under these conditionings,  $(X_t)_{t=\beta_i}^{\kappa_i}$ is just simple random walk conditioned to exit $\{\gamma_i\} \times H$ towards $\{\gamma_{i+1}\} \times H$. The effect of this conditioning can be understood via the Bayesian factors of~(\ref{e.Bayes}): from any vertex of  $\{\gamma_i\} \times H$, there is one edge towards the bag $\{\gamma_{i+1}\} \times H$, there are $2k-1$ edges towards other bags (all are forbidden by the conditioning), and $k^9-1$ edges inside $\{\gamma_i\} \times H$, hence the Bayesian factor of the edge towards $\{\gamma_{i+1}\} \times H$ is $2k$ times the Bayesian factor of each edge inside $\{\gamma_i\} \times H$. That is, $(X_t)_{t=\beta_i}^{\kappa_i}$ is just a simple random walk in the complete graph $\{\gamma_i\} \times H$, exiting after $\kappa_i-\beta_i  \stackrel{d}{=} \mathsf{Geom}\big(\frac{2k}{k^9-1+2k}\big)$ steps, independently of the actual steps taken inside $\{\gamma_i\} \times H$.

Now note that the process $L_t := \big| \LERW_t \cap (\{\gamma_i\} \times H) \big|$ for $t=\beta_i,\beta_i+1,\dots,\kappa_i-1$ is equal in distribution to the process of part (3) of Lemma~\ref{l.mixing}, started from the value at $t=\beta_i$ that is given by the conditioning. No matter what this starting value is, if $\kappa_i-\beta_i \geq k^5$, then the distribution of $L_{\kappa_i}$ is very close to the stationary distribution $\pi$:
\begin{align*}
\Pb{L_{\kappa_i} < k^3 \md  \calf_{i-1},\ \beta_r<\infty } &<  \Pb{\kappa_i-\beta_i < k^5} +   \pi(\{1,\dots,k^3\}) + d_\mathrm{TV}(\mu_{k^5},\pi) \\
&<  \PB{ \mathsf{Geom}\big(\frac{2k}{k^9-1+2k}\big) < k^5 } + \frac{k^6}{2k^9} + \exp(-ck^5/k^{4.5})\\
&< O(k^{-3}),
\end{align*}
where the second term was bounded using Part~(1) of Lemma~\ref{l.mixing}. Thus, we have proved~(\ref{e.large}).

Now, conditionally on $\beta_r<\infty$, checking the bags $\{\gamma_i\} \times H$ one-by-one for $i=0,1,\dots,r-1$, whether they satisfy $\big|\LERW_{\kappa_i} \cap (\{\gamma_i\} \times H)\big| \geq k^3$, the bound~(\ref{e.large}) tells us that the number of bags that do not satisfy this is stochastically dominated by a $\mathsf{Binom}\big(r,Ak^{-3})\big)$ variable, for some $A<\infty$. Now, a standard exponential Markov inequality tells us that
$$
\PB{ \mathsf{Binom}\big(r,Ak^{-3}\big) \ge r/2 } < C\exp\left(- r \, \frac{1}{2}\log k^3\right) = C k^{-\frac{3}{2}r},
$$
and therefore, if we call a bag {\it good} if $\big|\LERW_{\kappa_i} \cap (\{\gamma_i\} \times H)\big| \geq k^3$, then 
\be\label{e.many}
\PB{\textrm{there are more than } r/2 \textrm{ good bags } \{\gamma_i\} \times H } > 1-C k^{-\frac{3}{2}r}.
\ee
An important point here is that the ``bad event'' has probability much less than $k^{-r}$.
\medskip

The second part of the proof of~(\ref{e.avoidall}) is to study what happens from $\beta_r$ until $\tau_{\fb}$. 

Observe that, for any $i=1,\dots,r-1$, independently of all other bags, the set $(X_t)_{t=\varphi_i}^{\psi_i}  \cap  (\{\gamma_i\} \times H)$ stochastically dominates the set of vertices visited by a simple random walk on the complete graph $\{\gamma_i\} \times H$, run for $\mathsf{Geom}\big(\frac{2}{k^9+1}\big)$ steps. This is by comparing our set with the vertices visited by a random walk on a truncated graph, where the edges from $\{\gamma_i\} \times H$ towards ``side-bags'', bags other than $\{\gamma_{i\pm1}\} \times H$, are simply deleted. The comparison is simply by noting that every time we leave towards a side-bag, we also have to come back to $\{\gamma_i\} \times H$, possibly to the same vertex where we left, or to a uniformly random different vertex. And the distribution of steps on the truncated graph can be computed easily from the Bayesian factors (\ref{e.Bayes}): the conditioning that we have to reach $\{\gamma_0\} \times H$ before $\{\gamma_r\} \times H$ makes the weight of the step towards $\{\gamma_{i-1}\} \times H$ larger than the weight of the step towards $\{\gamma_{i+1}\} \times H$, but the sum of the two weights is still twice the weight of any of the steps within $\{\gamma_{i}\} \times H$, hence we will leave the bag $\{\gamma_{i}\} \times H$ on the truncated graph after $\mathsf{Geom}\big(\frac{2}{k^9+1}\big)$ steps.

So, every time the walk from $\beta_r$ until $\tau_{\fb}$ enters a good bag, the probability of $\LERW_{\kappa_i} \cap (X_t)_{t=\varphi_i}^{\psi_i} = \emptyset$ is bounded above (independently of what happened before) by the probability that $\mathsf{Geom}\big(\frac{2}{k^9+1}\big)$ independent steps avoid the at least $k^3$ vertices of $\LERW_{\kappa_i} \cap (\{\gamma_i\} \times H)$. This probability is at most $O(k^{-3})$. Therefore, if we denote the event of~(\ref{e.many}) by $\mathsf{Many}$, then
\begin{align*}
\PB{\bigcap_{i=1}^{r-1} \Avoid_i} &\leq \Pb{\mathsf{Many}^c} + \PB{\bigcap_{i=1}^{r-1} \Avoid_i \md \mathsf{Many}} \\
&\leq C k^{-\frac{3}{2}r} + O(k^{-3})^{r/2} = O(k^{-\frac{3}{2}r}).
\end{align*}
For $k$ large enough, this proves~(\ref{e.avoidall}), and we are done.
\qed

\section{Disco lights}\label{s.nonunimod} 

\begin{figure}[htbp]
\vglue 0.3 cm
\SetLabels
(0.6*1) $o=\gamma_0$\\
(0.28*0.77) $\gamma_1$\\
(0.13*0.52) $\gamma_2$\\
(-0.01*0.31) $z=\gamma_3$\\
(0.45*0.8) $\Delta_{0,0}$\\
(0.82*0.8) $\Delta_{0,1}$\\
(0.76*0.96) $\Delta_{0,2}$\\
\endSetLabels
\centerline{
\AffixLabels{
\includegraphics[width=0.95\textwidth]{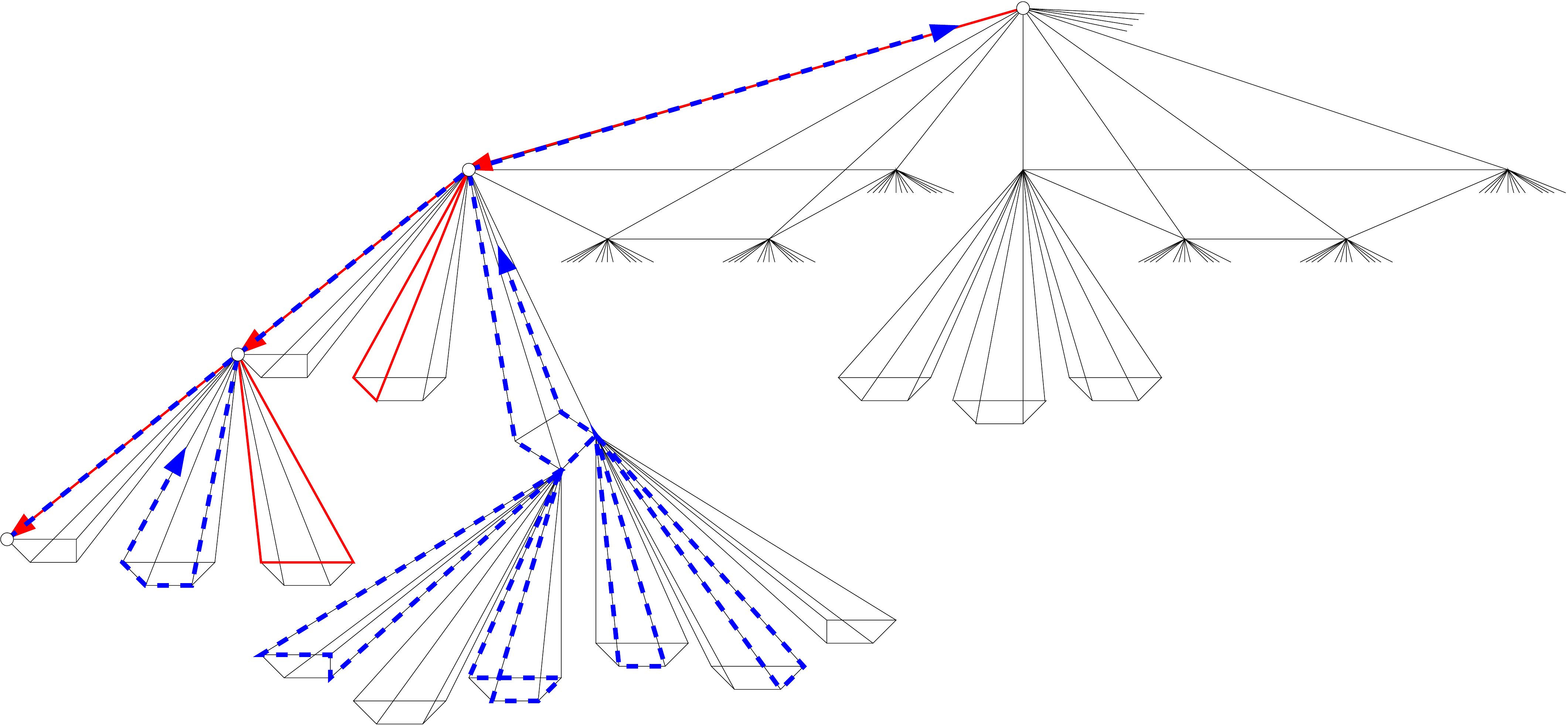}
}}
\caption{A random walk excursion in the pyramid graph $\Py^3$ that satisfies the good event $\cB_z$. The red solid parts are on the way to $z$, the blue dashed ones are on the way back.}
\label{f.diamonds}
\end{figure}
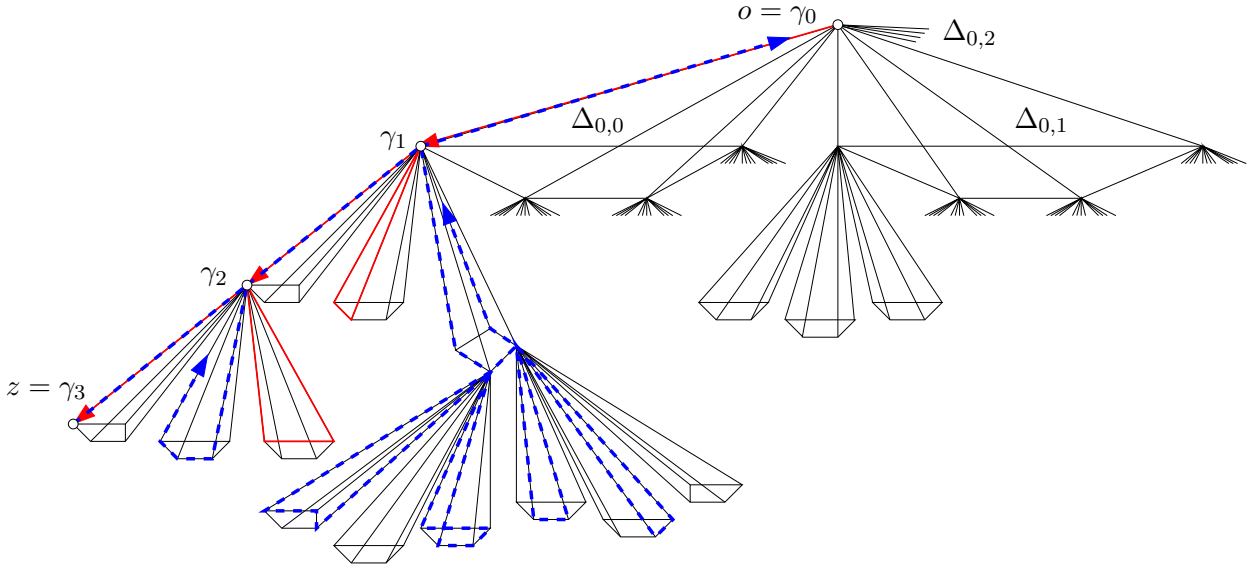

In place of $\T^k$, we will consider a transitive non-unimodular graph which we call the $k$-ary pyramid graph $\Py^k$. One pyramid is just a cycle $C_4$ with an extra vertex, the apex, connected to every vertex of this cycle. Now we take the tree $\T^{4k+1}$, and orient all of its edges towards a fixed end of the tree. For each vertex, divide the $4k$ incoming edges into 4-tuples, and connect the tails of the edges with a $C_4$. The resulting graph is $\Py^k$, which can also be considered as glued together from pyramids. See Figure~\ref{f.diamonds}. Then our example will be $G=\Py^k \times H$ for a large enough finite transitive graph $H$. This is obviously a transitive non-unimodular graph: if $\Gamma$ is the full automorphism group of $G$, and $(x,y)$ is an edge  of $\Py^k$ where $x$ is the apex of a pyramid and $y$ is in the base, then $|\Gamma_{(x,h)}(y,h)|=4k$, while $|\Gamma_{(y,h)}(x,h)|=1$, for any $h\in H$.

\bpr[Disconnected nonunimodular $\FUSF$]\label{p.discoPy}
For any $d\ge 2$, if $k$ is large enough, and $H$ is a connected finite $d$-regular transitive graph on at least $k^{5/2}$ vertices, then $\FUSF$ on $\Py^k\times H$ a.s.~has infinitely many components.
\epr

\proof Fix a geodesic ray $\gamma=\gamma(z)=\{o=\gamma_0,\gamma_1,\dots,\gamma_n=z\}$ from $o$ to $z \in S_n$ as before, let $\Delta_{i,0}$ be the pyramid containing both $\gamma_{i}$ and $\gamma_{i+1}$, and let $\Delta_{i,j}$, $j=1,\dots,k-1$, be the other pyramids with their apices at $\gamma_i$. Letting $(Y_t)_{t\ge 0}$ be simple random walk on $\Py^k$, started at $o$, the version of the event $\cB_z$ from~(\ref{e.avoid}) will be as follows. 
\be\label{e.ALPy}
\begin{aligned}
\cA_z &:= \Big\{ \textrm{the edge }(\gamma_{i-1},\gamma_{i})\textrm{ is crossed exactly twice by }(Y_t)_{t=0}^{\tau^+_o},\\
&\qquad\ \textrm{ and no other edge of }\Delta_{i-1,0}\textrm{ is crossed, for all } i=1,2,\dots,n\Big\},\\
\cL_z &:= \Big\{ \big|\big\{t\in \{1,\dots,\tau_z\} : Y_t=\gamma_i \big\} \big| \leq k/2+1 \textrm{, for all } i=1,2,\dots,n\Big\}.
\end{aligned}
\ee
%
Furthermore, let $E_i$ and $F_i$ be the set of pyramid bases $\Delta_{i,j}\setminus \{\gamma_i\}$ visited by $(Y_t)_{t=0}^{\tau_z}$ and by $(Y_t)_{\tau_z}^{\tau^+_o}$, respectively, among $j=1,\dots,k-1$, and  
define the event
\be\label{e.avoidPy}
\cB_z := \cA_z \cap \cL_z \cap \big\{E_i \cap F_i = \emptyset \textrm{ for all } i=1,2,\dots,n \big\}.
\ee

With these definitions, the proofs of Sections~\ref{s.tree},~\ref{s.product} go through almost verbatim, with three minor differences. The first one is that all the edges of the pyramids $\Delta_{i,0}$ except for $(\gamma_i,\gamma_{i+1})$ are forbidden for the random walk, which changes some probabilities by a uniform constant factor, on each level $i$. The second difference is that the graph of pyramids has now the tree structure (and thus, e.g., the self-avoidance condition $E_i\cap F_i=\emptyset$ is defined via pyramids), but the walk still chooses edges, not pyramids. Again, this can change  probabilities only by uniform constant factors. For instance, in place of~(\ref{e.There}) and~(\ref{e.pk}), we have
\be\label{e.TherePy}
\Pb{\mathsf{There}_i, \textrm{ after }j\textrm{ excursions}} = \left(1-\frac5k\right)^j \frac{1}{4k+1}
\ee
and
\be\label{e.pkPy}
p_k:= \Pb{\mathsf{There}_i,\ \mathsf{Back}_i, \textrm{ and } E_i \cap F_i = \emptyset } 
\ge \sum_{j=0}^{\lfloor k/2\rfloor} \left(1-\frac5k\right)^j\frac{1}{4k+1} \frac{1}{4j+2} 
\asymp \frac{\log k}{k}\,,
\ee
and everything works just as before.

The last minor difference is in the direct proof of having infinitely many trees almost surely, at the end of Section~\ref{s.product}. In the non-unimodular $\Py^k$, not all rays $o_1,o_2,\dots$ are the same; pick one tending to the distinguished end. Then, it is not obvious that the $\Py^k$-coordinate of the simple random walk,
viewed only at the times when it moves on this ray, is a one-dimensional {\it symmetric} walk. But it is, since the effective conductance between the cutpoints $o_i$ and $o_{i+1}$ is obviously equal to the effective conductance between $o_i$ and $o_{i-1}$, and hence, by a standard correspondence between hitting probabilities and electric networks \cite[Chapter 2]{LPbook}, we have $\Pso{o_i}{\tau_{o_{i-1}} < \tau_{o_{i+1}} }=1/2$.
\qed

Before proving Theorem~\ref{t.light}, we give our second proof of having infinitely many trees in the $\FUSF$ in the context of tree-like graphs that may also be nonunimodular. The next claim does not include any randomness, and for unimodular transitive graphs it is a tautology. 

\begin{proposition}[Infinite weight sums]\label{p.fewpieces}
Let $G$ be a transitive graph, $x\in V(G)$ fixed, and $S$ a finite set of vertices such that every component of $G\setminus S$ is infinite. Denote by $\Gamma$ the automorphism group of $G$, and by $\Gamma_y$ the stabilizer of a vertex $y$.
Then $\sum_{y\in C} \frac{|\Gamma_y x|}{|\Gamma_x y|}$ is infinite for every component $C$ of $G\setminus S$.
\end{proposition}

\proof
Let $M>1$ be the maximum of $\frac{|\Gamma_x y|}{|\Gamma_y x|}$ attained over neighbors $y$ of $x$. For a vertex $v$, let $N_1(v)$ be the set of all neighbors $y$ of $v$ with $\frac{|\Gamma_v y|}{|\Gamma_y v|}=M$. 
Note that if $y\in N_1(v)$ then $\gamma(y)\in  N_1(\gamma(v))$ for every $\gamma\in\Gamma$.
Define recursively $N_i(x)=\bigcup_{y\in N_{i-1}(x)} N_{1} (y)$ as $i=2,3,\ldots$. We have $|N_1 (x)|\geq M$, because $\Gamma_x y\subset N_1 (x)$ and $|\Gamma_y x|\geq 1$. 
We will see next $|N_i (x)|\geq M^i$. 

Choose an arbitrary $y\in N_i(x)$. Pick an arbitrary $\gamma\in\Gamma_x$, and fix a sequence $x_j\in N_1(x_{j-1})$ for $j=1,\ldots, i$ with $x_i=y$ and $x_0:=x$. We will prove by induction that $\gamma(x_j)\in N_j(x)$. For $j=1$ we have seen this. Then,
$x_{j}\in N_1(x_{j-1})$ implies $\gamma (x_{j})\in N_1(\gamma (x_{j-1}))$, and we know $N_1(\gamma (x_{j-1}))\subset N_1(N_{j-1}(x))=N_j(x)$ from the induction hypothesis. This completes the proof of $\Gamma_x y\subset N_i (x)$.
Apply the cocycle identity $\frac{|\Gamma_u v|}{|\Gamma_v u|}\frac{|\Gamma_v w|}{|\Gamma_w v|}=\frac{|\Gamma_u w|}{|\Gamma_w u|}$ to obtain $|N_i (x)|\geq |\Gamma_x y|=  {|\Gamma_y x|} \frac{|\Gamma_{x_0} x_1|}{|\Gamma_{x_1} x_0|}\cdots \frac{|\Gamma_{x_{i-1}} x_i|}{|\Gamma_{x_i} x_{i-1}|}\geq M^i$, as claimed.



For $v\in V(G)$, let $N_\infty (v):=\bigcup_{i=1}^\infty N_i (v)$, and let $m:=\min \Big\{\frac{|\Gamma_s x|}{|\Gamma_x s|} : s\in S\Big\}$. If there exists a $v\in C$ such that $\frac{|\Gamma_v x|}{|\Gamma_x v|}<m$, then $N_\infty (v)\cap S=\emptyset$ (using the simple observation that $\frac{|\Gamma_{{v'}} x|}{|\Gamma_x v'|}=\frac{|\Gamma_{v'} v|}{|\Gamma_v v'|}\frac{|\Gamma_{v} x|}{|\Gamma_x v|}=M^{-i}\frac{|\Gamma_{v} x|}{|\Gamma_x v|}< m$ for every $v'\in N_\infty (v)$, with some $i\geq 0$),
and hence 
$N_\infty (v)\subset C$. 
Then 
\begin{align*}
\sum_{y\in C} \frac{|\Gamma_y x|}{|\Gamma_x y|} \geq \sum_{y\in N_\infty (v)} \frac{|\Gamma_y x|}{|\Gamma_x y|} &=
\frac{|\Gamma_v x|}{|\Gamma_x v|}
\sum_{y\in N_\infty (v)} \frac{|\Gamma_y v|}{|\Gamma_v y|}\\
&=\frac{|\Gamma_v x|}{|\Gamma_x v|}
\sum_{i=1}^\infty\sum_{y\in N_i (v)} M^{-i}=\frac{|\Gamma_v x|}{|\Gamma_x v|}
\sum_{i=1}^\infty |N_i(v)|M^{-i}\\
&\geq \frac{|\Gamma_v x|}{|\Gamma_x v|} \sum_{i=1}^\infty 1.
\end{align*}
If there were no such $v$, then one would have an infinite sum of numbers at least $m$ in $\sum_{y\in C} \frac{|\Gamma_x y|}{|\Gamma_y x|}$, leading again to the conclusion that the sum is infinite.
\qed

Say that $\calt$ is a {\it tree-like decomposition} of a graph $G$ if $\calt$ is a random partition of $V(G)$ into finite connected sets, called {\it bags}, together with a tree on the bags such that any edge of $G$ goes between points of adjacent bags or within the same bag.  
We call a tree-like decomposition {\it invariant} if its distribution is preserved by the automorphisms of $G$. 

A random spanning forest $\calf$ of a graph $G$ was defined to be {\it weak insertion tolerant} in \cite{indist} if it satisfies the following property. Fix $r>0$ and vertices $x$ and $y$ of $G$ connected by an edge $e$. Let $\cD$ be the event that $x$ and $y$ are in different components of $\calf$. Then one can map every configuration $\omega\in \cD$ to a new configuration  $\omega \cup \{e\}\setminus \{f\}$, where $f$ is either the empty set or it is an edge of $\calf$ at distance at least $r$ from $x$, and it is determined by $\omega$ in a measurable way. The mapping just defined is measurable, and it takes events of positive probability (contained in $\cD$) to events of positive probability (contained in $\cD^c$). See \cite{indist} for a more thorough definition and the proof that the $\FUSF$ and the $\FMSF$ are weak insertion tolerant.

\begin{proposition}[1 or $\infty$ law]\label{p.inftymany}
Let $G$ be an infinite transitive graph that has an invariant random tree-like decomposition.
Let $\calf$ be an invariant random forest of $G$ with only infinite components, and suppose that it is weak insertion tolerant. Then $\calf$ has either 1 or infinitely many components almost surely.
\end{proposition}

\proof Proving by contradiction, suppose that $\calf$ has $m$ clusters with $1<m<\infty$. Let $\calt$ be an invariant random tree-like decomposition of $G$. We claim that with positive probability there exists a bag $B$ and a cluster $\cC$ of $\calf$ such that $B \cap \cC=\emptyset$. To see this, pick
a finite $B_0\subset G$ such that $B_0$ is a bag in $\calt$ with positive probability. Condition on this event, and on that some fixed and adjacent $x$ and $y$ in $B_0$ are in different $\calf$-clusters. (We may assume that $B_0$ was chosen so that this event has positive probability.)
Let $r>0$ be such that with probability at least $1/2$ any two points of $B_0$ that are in the same $\calf$-cluster have distance less than $r$ in $\calf$. Applying weak insertion tolerance, insert the edge $\{x,y\}$ with the possible removal of an edge at distance at least $r$ from $x$. Repeating this as many times as necessary (at most $|B_0|-1$ many times) for some other adjacent pairs in $B_0$, we arrive at an event of positive probability where all vertices of $B_0$ are in the same $\calf$-cluster. But then every other $\calf$-cluster has to be fully contained in one of the components of $G \setminus B_0$, and hence it cannot intersect the bags in the other components.  
Thus we have found some $B$ and $\cC$ as claimed.


For any $v\in V(G)$, if the bag $B_v$ of $v$ does not intersect some $\calf$-cluster $\cC$, then removing from $\calt$ all the bags that intersect $\cC$, 
we get some components, exactly one of which contains $B_v$. There is a unique $\calt$-edge from a unique bag $B_{v,\cC}^*$ of this component, to a bag that intersects $\cC$. (Possibly $B_{v,\cC}^*=B_v$.) Now, define the following mass transport: for $v,w \in G$,
$$
f(v,w,\calf) := \sum_{\cC} f_{\cC}(v,w,\calf,\calt),\qquad
f_{\cC}(v,w,\calf,\calt) :=
\begin{cases}
{1}/{|B_{v,\cC}^*|} & \text{if }w \in B_{v,\cC}^*\text{ for cluster }\cC\text{ of }\calf,\\
0 & \text{otherwise}.
\end{cases}
$$
We will use the {\it Tilted Mass Transport Principle} for invariant percolations on not necessarily unimodular transitive graphs \cite[Corollary 8.8]{LPbook}: if $\Gamma$ is the automorphism group of $G$, then
\be\label{e.MTP1}
\sum_{z\in V(G)} \E f(v,z,\calf,\calt) = \sum_{y\in V(G)} \E f(y,w,\calf,\calt) \, \frac{|\Gamma_y w|}{|\Gamma_w y|}.
\ee
The left hand side, which is the expected mass sent out, is clearly at most $m$. 

To estimate the right hand side, condition on the event $\cA$ that a fixed set $B$ is a bag of $\calt$ and it does not intersect some $\cC$ but is adjacent to a bag that intersects $\cC$. By the first paragraph in the proof, this has a positive probability. Fix a vertex $x$ in $B$.
Vertices $y$ in all but one component of $\calt \setminus B$ have the property that $B_{y,\cC}^*=B$, hence they all send mass $1/|B|$ to $x$. Furthermore, every infinite component of $\calt\setminus B$ contains an infinite component of $G\setminus B$, hence the right hand side of \eqref{e.MTP1} can be bounded from below by $\frac{\P(\cA)}{|B|}\E \sum_{y}  \frac{|\Gamma_y x|}{|\Gamma_x y|}$, where $y$ is running over the vertices in some infinite component of $G \setminus B$. (Which infinite component, that may depend on $\calt$.) By Proposition~\ref{p.fewpieces}, this sum is always infinite, leading to a contradiction to \eqref{e.MTP1}.
\qed

Once that the $\FUSF$ in $\Py^k\times H$ is disconnected with positive probability, Proposition~\ref{p.inftymany} gives that it has infinitely many components a.s.~by the ergodicity of the $\FUSF$ \cite[Section 10.4]{LPbook}.

\proofof{Theorem~\ref{t.light}} We already know that there are infinitely many trees in the $\FUSF$, and want to show that some of them are light. The fixed end of $\Py^k$ yields a natural projection $\pi: \Py^k\times H \longtwoheadrightarrow \Z$, where all the preimages $x\in\pi^{-1}(m)$ for a fixed $m\in \Z$ have the same Haar weight $\mu(\Gamma_x)=(4k)^{-m}$.

If all the infinitely many clusters $\cC_i$ in the $\FUSF$ were reaching infinitely high up in Figure~\ref{f.diamonds}, i.e., if $\inf\pi(\cC_i)=-\infty$, then, for any two $\cC_i, \cC_j$ there would exist some $\gamma_{m(i,j)} \in \pi^{-1}(m(i,j))$ such that the infinite geodesic ray $\gamma_{m(i,j)},\gamma_{m(i,j)-1},\dots$ in $\Py^k$, converging to the fixed end at $-\infty$, has the property that both clusters intersect each $\gamma_{m(i,j)-t}\times H$, for $t=0,1,\dots$. However, if we take enough clusters $\{\cC_i\}_{i\in I}$ such that ${|I|\choose 2} > |H|$, and let $m:=\min\{ m(i,j) : i,j \in I\}$, then $\gamma_m$ is already defined for each pair, and it is actually the same vertex of $\Py^k$, so we would need to have ${|I|\choose 2}$ disjoint clusters intersecting $\gamma_m \times H$, a contradiction.

Thus, all but finitely many clusters $\cC_i$ of $\FUSF$ have a smallest label $\min \pi(\cC_i) > -\infty$. Let $M(\cC_i)\subset G$ be the set of vertices achieving this minimal label; we set $M(\cC_i)=\emptyset$ if $\inf \pi(\cC_i) = -\infty$. Note that $|M(\cC_i)|\leq H$ almost surely. Now, define the following mass transport: for $x,y \in G$,
$$
f(x,y,\FUSF) :=
\begin{cases}
1 & \text{if }x,y\text{ are in the same component }\cC_i\text{ of }\FUSF, \text{ and }y\in M(\cC_i),\\
0 & \text{otherwise}.
\end{cases}
$$
We again use the {\it Tilted Mass Transport Principle} from \cite[Corollary 8.8]{LPbook}:
\be\label{e.MTP}
\sum_{y\in V(G)} \E f(x,y,\FUSF) = \sum_{y\in V(G)} \E f(y,x,\FUSF) \, \frac{\mu(\Gamma_y)}{\mu(\Gamma_x)}.
\ee
The left hand side is at most $|H|$. The right hand side, if $x\in M(\cC_i)$ for some cluster $\cC_i$, is 
$$
(4k)^{\min \pi(\cC_i)} \sum_{y\in \cC_i} \mu(\Gamma_y).
$$
By~(\ref{e.MTP}), this is finite, hence, whenever $\min \pi(\cC_i) > -\infty$, the cluster $\cC_i$ is light.\qed

\section{Tell me why}\label{s.open}

%
%

The first natural question is how general the phenomena of Theorems~\ref{t.disco} and~\ref{t.genset} really are: 

\begin{problem}
If $\Gamma$ is a finitely generated treeable group with $\WUSF \not= \FUSF$, does it always have two generating sets such that the $\FUSF$ is disconnected in the first Cayley graph, while it is connected in the second?
\end{problem}

Tom Hutchcroft has suggested that $\Z^5 * \Z^5$ might be a counterexample, where the disconnectedness of the $\FUSF$ cannot be destroyed. However, a proof  does not seem to be a trivial matter to us.

An affirmative answer for the connected case would of course imply $\beta^{(2)}_1(\Gamma)=\mathsf{cost}(\Gamma)-1$ for treeable groups, which is actually known to hold by Gaboriau's results \cite[Corollaries 3.23 and 3.16]{Gaboriau}.

For our specific Cayley graphs of Theorem~\ref{t.disco}, the following two problems remain open. Of course, a negative answer to Problem~\ref{p.touch} would be pointing towards a positive answer to Problem~\ref{p.add}.

\begin{problem}\label{p.touch} 
Is it true that if the $\FUSF$ in some $\T^k\times H$ is disconnected, then any two components touch each other only at finitely many places? Are there at least special choices for $k$ and $H$ for which this happens?
\end{problem}

\begin{problem}\label{p.add}
For the $\FUSF$ in any $\T^k\times H$, is the union of the $\FUSF$ with an independent $\mathsf{Bernoulli}(\eps)$ bond percolation connected, for any $\eps>0$? If not, is there any invariant way to make the $\FUSF$ connected by adding an arbitrarily small density edge percolation?
\end{problem}

As we already defined in the introduction, for any finite graph $H$ we let
$$
\disco(H):=\min\big\{k : \FUSF(\T^k \times H) \textrm{ is disconnected}\big\} \in \{3,4,\dots,\infty\}.
$$
We know that $\disco(P_2)=\infty$ from \cite{Pengfei}, while Theorem~\ref{t.disco} implies that if $\ell$ is large enough, then the cycle $C_\ell$ of length $\ell$ has $\disco(C_\ell) < \infty$. 

\begin{problem} 
What is the smallest $\ell$ for which $\disco(C_\ell) < \infty$? In particular, what is $\disco(C_3)$?
\end{problem}

\begin{problem} 
Are there infinitely many finite graphs $H$ with $\disco(H)=\infty$?
\end{problem}

The two choices for $H$ in Theorem~\ref{t.genset} inspire the following question:

\begin{problem} 
If $H$ and $H'$ are two finite connected graphs on the same vertex set and $E(H)\subset E(H')$, do we always have $\disco(H) \leq \disco(H')$?
\end{problem}

Far from this monotonicity, but a similar type of result was shown in \cite{ABIT}: for any $\T^k\times H$, if we increase the weights on the edges of $H$ enough, then the $\FUSF$ becomes connected.

Regarding the generality of Lemma~\ref{l.mixing}, we have not found the following question addressed in the literature. One piece of motivation is \cite{PeRe}.

\begin{problem} 
Is it true on any connected transitive graph on $n$ vertices that the typical size of the stationary loop-erased version of a simple random walk trajectory is $\Omega(\sqrt{n})$?
\end{problem}

Now, given how the proof of Theorem~\ref{t.disco} used Lemma~\ref{l.heatkernel}, and how the proof of Theorem~\ref{t.genset} used Lemma~\ref{l.mixing}, one may guess that if $H$ has better mixing properties, then disconnection becomes easier:

\begin{problem} 
If $H$ and $H'$ are two connected transitive $d$-regular finite graphs on the same vertex set, with $H'$ having a spectral gap larger than $H$, does it follow that $\disco(H) \geq \disco(H')$?
\end{problem}

The next natural player appearing on the floor is $\disco^*$, a parameter that is dual, in some sense, to $\disco$. Let us fix a natural sequence of finite graphs $\cH=(H_n)_{n\ge 1}$; as the simplest case, think of the cycles $H_n=C_n$. Then let
$$
\disco^*_\cH(k):= \min\big\{n : \FUSF(\T^k \times H_n) \textrm{ is disconnected}\big\}.
$$

\begin{problem} 
Consider the sequence of cycles $\cC=(C_n)_{n\ge 1}$. Is it the case that $\disco^*_\cC(3)<\infty$?
\end{problem}

\begin{problem} 
How about monotonicity in $k$? That is, if $\FUSF(\T^k \times H)$ is disconnected, then is $\FUSF(\T^{k+1} \times H)$ also disconnected?
\end{problem}

One can also define a continuous version of the graph parameter $\disco$. Recall from \cite{urn} or \cite[Chapter 14]{PGG} what unimodular random rooted graphs are.
\begin{align*}
\widetilde{\disco}(H) := \inf &\big\{\kappa : \FUSF(\cT \times H) \textrm{ is disconnected with positive probability}, \\
&(\cT,o)\textrm{ is an infinite unimodular random rooted tree with }\E \deg_{\cT}(o)=\kappa\big\} \in [2,\infty].
\end{align*}

\begin{problem} 
Find $\widetilde{\disco}(C_\ell)$.
\end{problem}

\begin{problem} 
Is there any finite graph $H$ with $\widetilde{\disco}(H) < \infty = \disco(H)$? 
\end{problem}

\begin{problem} 
Is there any finite graph $H$ with $\widetilde{\disco}(H) = 2$? (Note that  if $\E \deg_{\cT}(o)=2$, then $\cT$ has at most two ends, hence $\cT\times H$ is recurrent, hence the $\FUSF$ is connected almost surely. However, this does not exclude the possibility of the infimum being 2.) Or perhaps $\widetilde{\disco}(H) < \infty$ implies $\widetilde{\disco}(H)=2$?  
\end{problem}


\section*{Tributes}

We are indebted to Russ Lyons and Pengfei Tang for comments and corrections on the manuscript. We also thank Tom Hutchcroft and P\'eter Mester for useful remarks, and Damien Gaboriau and Asaf Nachmias for some references. Our work was supported by the ERC Consolidator Grant 772466 ``NOISE''. The second author was partially supported by Icelandic Research Fund Grant 185233-051.

\vskip 0.8 cm

\noindent {\bf G\'abor Pete}\\
Alfr\'ed R\'enyi Institute of Mathematics, Budapest,\\
and Institute of Mathematics, Budapest University of Technology and Economics\\
\url{http://www.math.bme.hu/~gabor}\\
{\tt gabor.pete [at] renyi.hu}
\vskip 0.4 cm

\noindent {\bf \'Ad\'am Tim\'ar}\\
Division of Mathematics, The Science Institute, University of Iceland\\
and
Alfr\'ed R\'enyi Institute of Mathematics\\
{\tt madaramit [at] gmail.com}

\end{document}